\documentstyle[12pt]{amsart}

\theoremstyle{plain}

\newtheorem{CorOwnNum}{Corollary}
\newtheorem{AssOwnNum}{Assumption}

\newtheorem{ThNoNum}{Theorem}

\newtheorem{ThNum}{Theorem}
\newtheorem{LemNum}[ThNum]{Lemma}
\newtheorem{PropNum}[ThNum]{Proposition}
\newtheorem{CorNum}[ThNum]{Corollary}

\theoremstyle{definition}

\newtheorem{NotNoNum}{Notation}
	 
\newtheorem{DefNoNum}{Definition}
	 
\newtheorem{ExamNoNum}{Example}
	 
\newtheorem{RemNoNum}{Remark}

\numberwithin{ThNum}{subsection}
\numberwithin{equation}{section}
\numberwithin{figure}{section}

\newcommand{\hyp}{{\bold H}}
\newcommand{\sphere}{{\bold S}}
\newcommand{\euc}{{\bold E}}
\newcommand{\integer}{{\bold Z}}
\newcommand{\rational}{{\bold Q}}
\newcommand{\real}{{\bold R}}
\newcommand{\complex}{{\bold C}}
\newcommand{\ball}{\operatorname{B}}
\newcommand{\Isom}{\operatorname{Isom}}
\newcommand{\area}{\operatorname{area}} 
\newcommand{\volume}{\operatorname{vol}}
\newcommand{\inj}{\operatorname{inj}} 
\newcommand{\length}{\operatorname{length}}
\newcommand{\radius}{\operatorname{radius}}

\newcommand{\Hom}{\operatorname{Hom}}

\newcommand{\PSL}{\operatorname{PSL}}
\newcommand{\SL}{\operatorname{SL}}

\newcommand{\arccosh}{\operatorname{arccosh}}

\begin{document}  

\title{%
	Deformations of Hyperbolic 3-Cone-Manifolds} 

\author{%
	Sadayoshi Kojima} 
\address{%
	Department of Mathematical and Computing Sciences \\
	Tokyo Institute of Technology \\
	Ohokayama, Meguro \\
	Tokyo 152 Japan} 
\email{%
	sadayosi@@is.titech.ac.jp}  

\begin{abstract}
	We show that any compact orientable  hyperbolic 
	3-cone-manifold with cone angles at most  $\pi$   
	can be continuously deformed 
	to a complete hyperbolic 
	manifold homeomorphic to the complement of the 
	singularity.  
	This together with the local rigidity by Hodgson 
	and Kerckhoff 
	implies the global rigidity for compact orientable  
	hyperbolic 3-cone-manifolds under the same angle 
	assumption. 
\end{abstract}

\maketitle

\setcounter{section}{-1}
\section{Introduction} 

A hyperbolic 3-cone-manifold is 
a riemannian 3-manifold of constant negative 
sectional curvature with cone-type 
singularity along simple closed geodesics 
(see \cite{HodgsonKerckhoff, Kerckhoff}).  
To each component of the singularity, 
associated is a cone angle.   
The cone angle is a positive real number.  
It possibly attains  $2\pi$.  
In this particular case, 
the singular set is not singular and 
simply a finite union of disjoint 
simple closed geodesics.  
The hyperbolic 3-cone-manifold is 
a generalization of the hyperbolic 3-orbifold 
with vertexless singularity.  

We are concerned with the deformations of 
a hyperbolic 3-cone-manifold with constant 
topological type.  
The hyperbolic Dehn filling 
theory by Thurston in  \cite{ThurstonNote},  
which describes deformations 
of a complete hyperbolic manifold 
in more wild setting, is 
a pioneering work of this subject.    
One important progress from a rather classical 
viewpoint was made recently by 
Hodgson and Kerckhoff  \cite{HodgsonKerckhoff}.  
They proved 
the local rigidity when 
cone angles are  $\leq 2\pi$.  
It corresponds to Weil \cite{Weil} and Garland 
\cite{Garland}  rigidity for hyperbolic manifolds 
of finite volume.  

In \cite{Kojima}, we showed as one of applications of 
the local rigidity that the underlying 
space of a nonsingular part of  
a compact hyperbolic 3-cone-manifold admits 
a complete hyperbolic structure of finite volume.  
We come back to this result later.  
The purpose of the present paper is to connect 
that complete structure with the original 
singular structure 
by a continuous family of cone-manifolds
with constant topological type 
under the assumption that initial cone angles are  $\leq \pi$.  
The main theorem is 

\begin{ThNoNum}\label{MainTh} 
Let  $C$  be a compact orientable hyperbolic 
3-cone-manifold and  $\varSigma$  a singular set which 
forms a link in  $C$.  
If the cone angles assigned to the components of  $\varSigma$  
all are at most  $\pi$,  then 
$C$  admits an angle decreasing continuous family of deformations  
to the complete hyperbolic manifold 
homeomorphic to the nonsingular part  $C-\varSigma$.  
\end{ThNoNum}

There are two immediate 
corollaries related to 
the representations 
in the group of orientation preserving 
isometries of the hyperbolic 3-space  $\hyp^3$, 
isomorphic to  $\PSL_2(\complex)$.  
For cone-manifolds, 
we have a holonomy representation of a nonsingular 
part, so that meridional loops of the singularity 
are mapped to elliptic elements. 
The representation so obtained could be 
neither faithful nor discrete.  
Nevertheless, the local rigidity in \cite{HodgsonKerckhoff} 
asserts 
that a neighborhood of this wild 
representation 
is parameterized up to conjugacy 
by Dehn filling coefficients  \cite{ThurstonNote}  
which are geometrically well understood.  

One corollary is about the global rigidity.  
Weil and Garland rigidity 
states that any nearby discrete faithful representations  
of a group in  $\PSL_2(\complex)$  with finite volume quotients  
are conjugate.  
In the case of cone-manifolds, 
the local rigidity says in particular that  
nearby representations of a holonomy representation of 
a nonsingular part with constant rotation angles 
for meridians are conjugate.  
Mostow \cite{Mostow} and  Prasad \cite{Prasad}  rigidity 
for a hyperbolic manifold of finite volume then asserts that
any discrete faithful representations 
in  $\PSL_2(\complex)$  are not only locally but 
globally conjugate each other.  
This global rigidity implies a geometric 
consequence that 
homeomorphic hyperbolic manifolds 
of finite volume are isometric.  
We state the global rigidity 
for hyperbolic 3-cone-manifolds  
rather in terms of this geometric terminologies.  

\begin{CorOwnNum}\label{GlobalRigidity}
Let  $C$  be a compact orientable hyperbolic 3-cone-manifold 
with singularity  $\varSigma$  where 
cone angles assigned to the components of 
$\varSigma$  all are  $\leq \pi$.  
If  $(C', \varSigma')$  is homeomorphic to  $(C, \varSigma)$  
so that the corresponding cone angles all are the same,  
then  $C'$  and  $C$  are isometric.   
\end{CorOwnNum}
 
The other corollary is 
about liftability of 
$\PSL_2(\complex)$-representations  
into  $\SL_2(\complex)$. 
The liftability has been discussed 
particularly for discrete subgroups in  $\PSL_2(\complex)$.  
As was pointed out in  \cite{Culler},  
the liftability depends only on the component of 
the space of representations.  
The main theorem will be proved by connecting 
the holonomy representation of a cone-manifold 
with that of 
a complete structure by a particular path 
in the space of representations.  
Since the holonomy representation of a complete 
hyperbolic manifold is known to lift in  \cite{ThurstonNote}, 
we have 

\begin{CorOwnNum}\label{HolonomyLifts}
The holonomy representation of a compact orientable 
hyperbolic 3-cone-manifold can be lifted to  
a  $\SL_2(\complex)$-representation if the cone angles 
assigned all are at most  $\pi$.  
\end{CorOwnNum}

It is quite unlikely that the angle assumptions in Theorem and 
Corollaries are necessary, though the argument 
we develop here uses its advantage.  
More progress should be expected.  

A fairly large part of the proof of the main theorem is due to 
Thurston's strategy for the 
geometrization of orbifolds \cite{symmetry,ThurstonBull} 
together with the local rigidity by Hodgson and 
Kerckhoff.  
The over all logic of Thurston's argument
and some of its details can be found in   
\cite{SOK}.  
We convey its minimal essentials for our purpose here, 
and hence 
the exposition will be reasonably self-contained.  

This paper is organized as follows.  
In the first section, we will review some 
basic facts about hyperbolic 3-cone-manifolds.  
Also we improve the results in  \cite{Kojima} 
from more Riemannian geometric viewpoints  
due to Steve Kerckhoff.  
The second section is to introduce two main 
machineries, the local rigidity and the pointed 
Hausdorff-Gromov topology.  
They are fundamental 
when cone angles are  $\leq 2\pi$. 
The third section is to establish a few tools 
to control the local geometry of 
cone-manifolds away from the singularity  
when cone angles are  $\leq \pi$.  
This section contains a technical but the most crucial part of 
the analysis.  
In section 4,  we establish a discrete setting of 
the problem.  
Then we study what happens when 
tubular neighborhoods of the singularity in the 
deformations are
uniformly thick in section 5, and 
when otherwise in section 6.  
In the final section,  we study what 
happens to continuous families and 
prove Theorem and Corollaries.  

The author would like to thank Steve Kerckhoff for 
careful attention to this work and for
showing him a quick idea to prove 
Theorem~\ref{ComplementHyperbolic} in full generality,  
Teruhiko Soma and the members of the Saturday Seminar at 
Tokyo Institute of Technology, especially Shigenori Matsumoto,   
for their invaluable suggestions, and the Centre Emile Borel 
of the Institut Henri Poincar\'e for their hospitality, 
where he finished the first version of this paper. 

In addition, the author would like to thank 
the referee for pointing out an earlier work by 
Qing Zhou \cite{Zhou} which discusses 
similar deformations.  
Zhou's argument together with the local 
rigidity now established by Hodgson and 
Kerchkhoff \cite{HodgsonKerckhoff} lead to 
another proof of the existence of an angle decreasing 
family in the main theorem provided 
that the initial cone angles all are at most  $2\pi/3$. 
Also. Zhou proved Theorem~\ref{ComplementHyperbolic} in 
slightly different manner.  


\section{Hyperbolic 3-Cone-Manifolds}

In this section, we assemble a few standard notions 
and notations which we use throughout this paper, 
improve the results in \cite{Kojima}  
and discuss an upper bound of the volume of 
hyperbolic 3-cone-manifolds with the same topological type.


\subsection{Preliminaries}

Let  $X$  be a metric space with a metric  $d$.  
An $R$-neighborhood of  $x \in X$  for  $R > 0$  
is the set of points in  $X$  
from which the distance to  $x$  is  $< R$,  
and denoted by 
\begin{equation*} 
 \ball_R(X, x) = \{ y \in X \, \vert \, d(y, x) < R \}.      
\end{equation*} 

If  $X$  is a riemannian manifold,  
the closure of  $\ball_R(X, x)$  is 
homeomorphic to a closed ball 
at least for sufficiently small  $R$.  
The injectivity radius of  $X$  at  $x \in X$  is 
the first supremum of such radii, and 
denoted by  $\inj_x X$.  
If  $X$  has a boundary  $\partial X$,  
we choose the supremum by furthermore requiring that 
$\ball_R(X, x)$  does not touch  $\partial X$.  
The injectivity radius of  $X$  is the infimum of 
injectivity radii of the points in  $X$.  
We denote it by 
\begin{equation*} 
\inj X \, (:= \inf \, \{ \inj_x X \, \vert \, x \in X \}).  
\end{equation*} 
The injectivity radius for manifolds with nonempty boundary 
by this definition would not be interesting since 
the points close to the boundary always make it 
vanishing.  

Let  $C$  be an 
orientable hyperbolic 3-cone-manifold of finite volume  
with compact singularity. 
The singular set  $\varSigma$    
is assumed to form a link  
\begin{equation*} 
\varSigma = \varSigma^1 \cup \cdots \cup \varSigma^n 
 \end{equation*} 
of  $n$ components.  
To each component  $\varSigma^j$  
of  $\varSigma$, associated 
is a cone angle  $\alpha^j \in (0, \infty)$.  
The angle set  $A$  of  $C$  is 
a vector   
\begin{equation*} 
A = (\alpha^1, \cdots, \alpha^n)
\end{equation*} 
of cone angles.  

$C$  carries a nonsingular but incomplete 
hyperbolic structure on the complement of the singularity
\begin{equation*} 
N = C - \varSigma.  
\end{equation*} 
$C$  itself inherits a metric induced from a Riemannian metric 
on  $N$.  
We assume that  $C$  is complete with respect to 
this metric.  
In particular, the metric completion of  $N$  is 
identical to  $C$.  
We have a  
developing map of  $N$  from its universal 
covering space  $\widetilde{N}$, 
\begin{equation*} 
{\cal D}_C: \widetilde{N} \to \hyp^3, 
\end{equation*} 
and a holonomy representation  
\begin{equation*} 
\rho_C : \Pi = \pi_1(N) \to \PSL_2(\complex). 
\end{equation*} 
They are called a developing map and a holonomy 
representation of a cone-manifold  $C$.  
A developing map is a local isometry, but never be injective.  
A holonomy representation is hardly discrete nor faithful.  

Let  $m_j, \; j = 1, 2, \cdots, n,$  be an oriented meridional 
loop for each component of  $\varSigma$.  
The image  $\rho_C(m_j)$  of a meridian  $m_j$  by the holonomy
representation is an elliptic element rotating  $\hyp^3$  by 
$\alpha^j $  about the axis,  though the rotation angle  
of  $\rho_C(m_j)$  makes sense 
only modulo  $2\pi$.  

The injectivity radius of  
$C$  at   $x \in N = C - \varSigma$  is 
to be the injectivity radius of  $N$  at  $x$  
and denoted by 
\begin{equation*}
\inj_x C \; (: = \inj_x N).  
\end{equation*}
The global injectivity radius of a 
cone-manifold  $C$  by this definition would not 
be interesting since 
if the singular set is 
nonempty, then the points close to the singularity 
always make it vanishing.  

\begin{DefNoNum} 
A {\it topological type} of a cone-manifold  $C$  is 
a homeomorphism type of 
a pair  $(C, \varSigma)$.  
We say  $C$  is {\it homeomorphic} to  $C_1$  for short 
if there is a homeomorphism between 
$(C, \varSigma)$  and  $(C_1, \varSigma_1)$.  
More strong relation is an {\it isomorphism type}.  
Two cone-manifolds are {\it isomorphic} if 
they share not only topological types but also 
cone angles, more precisely if there is 
a homeomorphism between  $(C, \varSigma)$  and  
$(C_1, \varSigma_1)$  so 
that the corresponding components of the singularity 
share the same cone angles.  
Such a homeomorphism is called an {\it isomorphism}.  
A self isomorphism is called an {\it automorphism} as usual. 
The strongest relation is an {\it isometry type} whose 
definition would be obvious.  
\end{DefNoNum} 

\begin{RemNoNum} 
The global rigidity is the claim that  
the isomorphism type and the isometry type are the same.  
\end{RemNoNum}


\subsection{Nonsingular Parts} 

The following theorem was proved in  \cite{Kojima} 
under an extra angle assumption 
using Hodgson-Kerckhoff's local rigidity.  
Here we present a quick argument due to Steve Kerckhoff  
which does not use the local rigidity and 
works without any angle assumption.  
As we mentioned in the introduction, 
Zhou also showed the following theorem in  \cite{Zhou}  
in slightly different manner.  

\begin{ThNum}\label{ComplementHyperbolic}
The underlying space of a nonsingular part  $N$  of 
an orientable hyperbolic 
3-cone-manifold  $C$  of finite volume 
carries a complete negatively curved metric.  
In particular it is 
homeomorphic to an interior of a 
compact irreducible atoroidal 3-manifold with toral 
boundary which admits no Seifert fibrations.  
Moreover, it admits a 
complete hyperbolic structure of finite volume.  
\end{ThNum}

\begin{pf} 
In the cylindrical coordinates around each component of 
the singularity  $\varSigma$, 
the metric has the form 
\begin{equation*} 
d \delta^2 + \sinh^2 \delta \, d \theta^2 + \cosh^2 \delta \, d \lambda^2  
\end{equation*} 
where  $\delta$  is the distance from the singularity, 
$\lambda$  is the distance along the singularity, 
$\theta$  is the angular measure around the 
singularity.  
Choose  $\varepsilon > 0$  small enough so that 
an $\varepsilon$-tubular neighborhood 
of  $\varSigma$ is a disjoint union of 
an $\varepsilon$-tubular neighborhood of each component 
of  $\varSigma$.  
Also choose 
monotone $C^{\infty}$-functions 
$\varphi(\delta)$  and  $\psi(\delta)$  in terms of  $\delta$  
so that 
\begin{equation*} 
	\varphi(\delta) = 
	\begin{cases} 
		\, 1 & \quad \text{if} \;\; \delta \geq \varepsilon, \\ 
		\, \text{O} (1/\delta) & \quad \text{if} \;\; \delta \to 0,
	\end{cases} 
\end{equation*}
and 
\begin{equation*} 
 	\psi(\delta) = 
	\begin{cases} 
	\; \cosh \delta  & \quad \text{if} \;\; \delta \geq \varepsilon, \\ 
	\; \text{O} (\delta) & \quad \text{if} \;\; \delta \to 0, 
	\end{cases} 
\end{equation*} 
where  $\text{O}( \; )$  is the Landau symbol, 
and modify the metric in an $\varepsilon$-tubular 
neighborhood of each component of  $\varSigma$  by 
\begin{equation*} 
\varphi^2(\delta) d\delta^2 + \sinh^2 \delta \, 
	d \theta^2 + \psi^2(\delta) d \lambda^2.   
\end{equation*} 
This gives a complete metric on the nonsingular part 
$N = C - \varSigma$  
since  $\varphi(\delta)$  diverges when  $\delta \to 0$.  

Let us compute 
the sectional curvature for this new metric.  
For notational convenience, 
we set  $\delta = x_1$, $\theta = x_2$  and $\lambda = x_3$.  
By a computation of the Christoffel symbols, 
we have 
the evaluation of the connection with respect to 
this basis,  
\begin{table}[h]
  \begin{center} 
	\begin{tabular}{c|ccc}
		$\nabla_{\partial/\partial x_i}(\partial/ \partial x_j)$ & 
		$1$ & $2$ & $3$ \\ 
	\hline	
		$1$ & $\varphi'/\varphi \cdot \partial/\partial x_1$ & 
		$\coth \delta \cdot \partial/ \partial x_2$ & 
		$\psi'/\psi \cdot \partial/\partial x_3$ \\ 
		$2$ & $\coth \delta \cdot \partial/\partial x_2$ & 
		$-\sinh \delta \cosh \delta/\varphi^2 \cdot 
			\partial/ \partial x_1$ & $0$ \\ 
		$3$ & $\psi'/\psi \cdot \partial/ \partial x_3$ & 
		$0$ & $-\psi \psi'/\varphi^2 \cdot \partial /\partial x_1$ 
	\end{tabular}
  \end{center} 
\end{table} 

\noindent
It is symmetric and 
hence the Riemannian curvature tensor 
\begin{equation*} 
R(X,Y,Z) = ([\nabla_X, \nabla_Y] - \nabla_{[X,Y]})(Z) 
\end{equation*} 
has been diagonalized, i.e., 
$<R(X,Y,Z),W> = 0$  unless  
$X \ne Y$  and either  
$(X,Y) = (Z,W)$  or  $(X,Y) = (W,Z)$.  
Thus the sectional curvatures 
\begin{equation*} 
K(X,Y) = - \frac{<R(X,Y,X),Y>}
	{\Vert X \Vert^2 \cdot \Vert Y \Vert^2 - <X,Y>^2} 
\end{equation*} 
all are convex combinations of the 
three sectional curvatures  
$K_{ij} = K(\partial/\partial x_i, \partial / \partial x_j)$  
with  $\{ i,j \} = \{ 1,2,3 \}$.  
By a direct computation, we have 
\begin{align*} 
	K_{12} 
	& = \frac{1}{\varphi^2} \left( \frac{\varphi' \coth \delta}{\varphi} -1 \right), \\
	K_{13} 
	& = - \frac{1}{\varphi^2} \left( \frac{\psi''}{\psi} 
		- \frac{\varphi' \psi'}{\varphi \psi} \right), \\ 
	K_{23} 
	& = - \frac{\psi' \coth \delta}{\varphi^2 \psi}.   
\end{align*} 
Then these three values stay negative away from zero  
for  $0 < \delta \leq \varepsilon$  
by choosing  $\varphi$  and  $\psi$  appropriately, 
and that  the volume change is still bounded.  
Hence the nonsingular part 
$N$  admits a complete 
negatively curved metric of finite volume.  

The topological properties in the 
statement all are now easily derived from the  
negativity of the new metric.  
Then by Thurston's uniformization theorem for 
Haken manifolds \cite{Morgan, ThurstonBull},  
$N$  admits a complete hyperbolic structure.  
\end{pf}

It is convenient to adopt the 
convention that we   
regard a complete hyperbolic 
manifold with preferred meridional elements 
as a cone-manifold.  

\begin{DefNoNum}
Let  $M$  be a complete orientable hyperbolic 3-cone-manifold 
of finite volume with prescribed meridional 
elements not only for the components of  $\varSigma$  
but for some cusps.  
We call  $M$  also a cone-manifold. 
The {\it topological type} of  $M$  is  
a pair of the Dehn filled resultant of  $M$  along 
prescribed meridional curves for cusps in question and 
the union of  $\varSigma$  and surgery cores.  
The cone angles which should be assigned to the 
components for cusps all are zero.  
\end{DefNoNum} 

\begin{NotNoNum}
If  $C$  is a compact orientable hyperbolic 
3-cone-manifold in the original sense, 
then the nonsingular part together with the 
complete structure obtained by 
Theorem~\ref{ComplementHyperbolic} is 
a cone-manifold in this new sense 
by assigning original 
meridional elements of  $C$.  
We let  $C_{comp}$  be 
this particular cone-manifold without singularity.  
The topological types of  $C$  and  $C_{comp}$  are 
the same. 
\end{NotNoNum}  

\begin{CorNum}\label{AutomorphismFinite}
The group of automorphisms of an 
orientable hyperbolic 3-cone-manifold   
of finite volume up to isotopy is finite.  
\end{CorNum}

\begin{pf} 
Here we duplicate the proof in  \cite{Kojima}.  
Any automorphism of  $C$  induces a homeomorphism 
of  $C_{comp}$  preserving meridians.  
By Mostow-Prasad rigidity, it is homotopic to 
an isometry of  $C_{comp}$.   
Moreover by Waldhausen's theorem \cite{Waldhausen}, 
this homotopy can be attained by an isotopy.  
The isotopy extends to an automorphism of  $C$.  
Hence the group of automorphisms of  $C$  up to 
isotopy is realized as a subgroup of the isometry 
group of  $C_{comp}$,  which is finite.  
\end{pf}

\begin{RemNoNum} 
The existence of the singularity in the proof is needed 
to apply Waldhausen's result, and in fact 
the argument does not cover the case without singularity.  
This finiteness up to isotopy for nonsingular case 
had been a difficult question and 
was settled very recently 
by Gabai, Meyerhoff and Thurston \cite{GMT} 
in full generality, 
though the finiteness up to homotopy is a consequence 
of Mostow rigidity.  
\end{RemNoNum}

\begin{CorNum}\label{EndToral} 
An orientable hyperbolic 3-cone-manifold 
of finite volume admits only finitely many 
toral cusps which are foliated by horotori. 
\end{CorNum}

\begin{pf} 
Since the end of  $C_{comp}$  consists of finitely many 
toral ends,  
the end of  $C$  consists of finitely many components 
homeomorphic to the torus time $\real$  at least topologically.  
To get a complete end, 
every element of a fundamental 
group of each end component must be mapped 
by the holonomy representation to a parabolic element 
with a common fixed point at the sphere 
at infinity  $\sphere^2_{\infty}$.  
A small neighborhood of each component is 
then foliated by the quotient of horospheres by 
the holonomy image. 
\end{pf}

We conclude this subsection by the following observation 
for loops in a hyperbolic 3-cone-manifold.  
Every loop in a hyperbolic manifold admits a length 
shortening free homotopy to either a closed geodesic 
(including a point)  
or an arbitrary short loop tending to the cusp.  
This will not be true for cone-manifolds in general, 
but we have its weaker version.  

\begin{LemNum}\label{LengthShorteningHomotopyExist} 
Any loop in the nonsingular part  $N$  of 
an orientable hyperbolic cone-manifold  $C$  of 
finite volume admits  
a length shortening homotopy to 
either a closed geodesic,  
an arbitrary short loop tending to the cusp or 
a piecewise geodesic loop hitting 
the singularity  $\varSigma$.  
\end{LemNum} 

\begin{pf} 
Let  $\ell$  be an arbitrary loop in the nonsingular 
part  $N$.  
Fix a point  $p$  on  $\ell$,  and we straight  $\ell$.  
Then it defines a length shortening homotopy 
to at least a piecewise geodesic loop hitting 
the singularity  $\varSigma$.  
If  $\varSigma$  does not obstruct the straightening, 
the length shortening homotopy reaches to a geodesic path  $\ell'$  
based at  $p$.  
We then homotop  $\ell'$  to a shorter geodesic path 
by sliding the reference point.  
This length shortening homotopy either terminates at 
some stage or tends to the cusp.  
When it terminates, 
either the homotopy hits  $\varSigma$  or 
it reaches to a closed geodesic.  
\end{pf}

\begin{RemNoNum} 
The proof is valid even when  $C$  is of infinite volume, 
or when the singularity is noncompact.   
\end{RemNoNum}


\subsection{Volumes}

The volume of a hyperbolic 
cone-manifold  $C$  can be expressed in terms of 
an integral over the canonical section of the associated 
flat  $\hyp^3$-bundle over  $N$.  
To say more precise, 
let  $\widetilde{N}$  be the universal covering space of  $N$.  
Each element  $\iota \in \Pi = \pi_1(N)$  acts on  
$\widetilde{N}$  by a deck transformation  
and on  $\hyp^3$  by 
$\rho_C(\iota) \in \rho_C(\Pi)$.  
Denote by  $E$  the fiber product  
$\widetilde{N} \times \hyp^3 / (\Pi \times \rho_C)$.  
It is a $\hyp^3$-bundle over  $N$  with a structure group 
in  $\PSL_2(\complex)$.  
The developing map  ${\cal D}_C: \widetilde{N} \to \hyp^3$  
defines a section  
$id \times {\cal D}_C : \widetilde{N} \to \widetilde{N} \times \hyp^3$  
which descends to a canonical one 
$s_0 : N \to E$  with respect to the hyperbolic structure on  $N$.  
The volume of  $C$  is then identified with the integral, 
\begin{equation}\label{SectionIntegral}
\volume C = \int_N s_0^* \, dv, 
\end{equation}
where  $dv$  is the volume form of  $\hyp^3$.  

\begin{LemNum} 
Suppose  $C$  is compact.    
If a section  $s_1 : N \to E$  agrees with 
the canonical one  $s_0$  on a small 
tubular neighborhood of  $\varSigma$, 
then 
\begin{equation*} 
\int_N s_1^* \, dv = \int_N s_0^* \, dv.  
\end{equation*}
In particular, the volume defined in \eqref{SectionIntegral} 
depends only on the behavior of the section 
near the singularity.  
\end{LemNum} 

\begin{pf}
Since  $\hyp^3$  is contractible, there is a homotopy 
$s_t$  of sections connecting  $s_0$  and  $s_1$.  
Let  $N_{\varepsilon}$  be a compact exterior 
of an $\varepsilon$-tubular neighborhood of   $\varSigma$.   
Choose  $\varepsilon$ small enough 
so that  $s_0 = s_1$  on that neighborhood of  $\varSigma$. 
Then by Stokes, 
\begin{align*} 
\int_N s_1^* \, dv - \int_N s_0^* \, dv 
	= & \int_{N_{\varepsilon}} s_1^* \, dv 
		-\int_{N_{\varepsilon}} s_0^* \, dv \\
	= & \int_{\partial N_{\varepsilon} \times [0, 1]} s_t^* \, dv \\
	= & \int_{N_{\varepsilon} \times [0, 1]} d \, (s_t^* \, dv). 
\end{align*} 
The last integral is zero.    
\end{pf}

\begin{PropNum}\label{Volume<Bound}
Given a compact orientable 
hyperbolic 3-cone-manifold  $C$  with  
the angle set  $A = (\alpha^1, \cdots, \alpha^n)$, 
there is a constant  $V_{max}$  
so that if a hyperbolic cone-manifold  $C_1$  
homeomorphic to  $C$  has an  
angle set  $A_1 = (\alpha^1_1, \cdots , \alpha^n_1)$  with  
$\alpha^j_1 \leq \alpha^j$  for all  
$1 \leq j \leq n$,  
then  $\volume C_1 \leq V_{max}$.  
\end{PropNum}

\begin{pf} 
Take a fine geodesic triangulation  $K$  of  $C$  so that 
$\varSigma$  is a subcomplex and 
a star neighborhood of  $\varSigma$  is 
a closed regular neighborhood.  
Then since  $\alpha^j_1 \leq \alpha^j$  for all  $j$, 
we may choose a homeomorphism  
$\varphi : (C, \varSigma) \to (C_1, \varSigma_1)$  so that 
any 3-simplex in a regular neighborhood of  $\varSigma$  
is mapped to an honest geodesic simplex in  $C_1$.  

Let  $\widetilde{K^{(0)}}$  be the preimage of the 0-skeleton 
$K^{(0)}$  in the universal cover  $\widetilde{N}_1$.  
Then the developing map  ${\cal D}_{C_1}$  defines the map  
$id \times {\cal D}_{C_1}\vert_{\widetilde{K^{(0)}}} : 
\widetilde{K^{(0)}} \to \widetilde{N}_1 \times \hyp^3$  
which extends to an equivariant continuous map of  $\widetilde{N}_1$  
by straightening.  
Since it is equivariant, it descends to a 
section  
$s_1 : N_1 \to E_1 = 
	\widetilde{N}_1 \times \hyp^3/(\Pi_1 \times \rho_{C_1})$.  
$s_1$  is identical to the canonical section 
$s_0 : N_1 \to E_1$  near  $\varSigma_1$.  
Hence by the previous lemma, we have  
\begin{equation*} 
\volume C_1 = \int_{N_1} s_1^* \, dv. 
\end{equation*}

The right hand side is a total sum of the signed volumes 
of 3-simplices appeared by straightening.  
However since the volume of a 3-simplex is uniformly 
bounded by a constant  
$v_3$  (= the volume of a regular ideal tetrahedra), 
this sum is bounded by  
$V_{max} = k v_3$,  
where $k$  is the number of 3-simplices in  $K$.  
The constant  $V_{max}$  obviously depends 
only on  $C$  and not on  $C_1$.   
\end{pf}

\begin{RemNoNum}
A number of tetrahedra needed to triangulate a 
neighborhood of the singularity  $\varSigma$  in  $C$  
depends in fact on the cone angles.  
However it is bounded by some constant 
depending only on the upper bound of  
the cone angles and not on any particular  $C$.  
\end{RemNoNum}


\section{Cone Angles $\leq 2 \pi$}

In this section, we review two machineries to study  
deformations of a hyperbolic 3-cone-manifold  $C$  
when cone angles are at most  $2 \pi$.  
One is the local rigidity by 
Hodgson and Kerckhoff \cite{HodgsonKerckhoff}  and 
the other is the pointed Hausdorff-Gromov topology 
studied in \cite{Gromov}.  


\subsection{Local Rigidity} 

Let us recall what the deformation of 
a hyperbolic 3-cone-manifold  $C$  is.  

\begin{DefNoNum} 
A deformation of a hyperbolic 3-cone-manifold  $C$  is 
a hyperbolic 3-cone-manifold  $C_1$  
together with a reference homeomorphism  
$\xi_1 : (C, \varSigma) \to (C_1, \varSigma_1)$.  
Two deformations  $(C_1, \xi_1)$  and 
$(C_2, \xi_2)$  of  $C$  are equivalent if 
there is an isometry  $\psi : C_1 \to C_2$  so that  
$\xi_2$  is isotopic to  $\psi \circ \xi_1$.  
A composition of a developing map  ${\cal D}_{C_1}$  with  
a lift  $\tilde{\xi}_1$  of  $\xi_1$  
defines a continuous map 
${\cal D}_{C_1} \circ \tilde{\xi}_1 : \widetilde{N} \to \hyp^3$.   
The developing map is well defined up to 
multiplication of isometries of  $\hyp^3$,  
and hence a deformation can be considered as a point on 
a quotient space of the mapping space 
${\cal M}(\widetilde{N}, \hyp^3)$  
by  $\PSL_2(\complex)$-action on the images.  
The quotient space carries a topology induced by 
a compact open topology on 
${\cal M}(\widetilde{N}, \hyp^3)$.  
The set of equivalence classes of deformations 
of  $C$  carries a further quotient topology by 
the action on the source of the group of lifts of 
automorphisms of  $C$  which are isotopic to the identity. 
We say  $(C_1, \xi_1)$  is a small deformation of  $(C, id)$  
if it is close to  $C$  in this topology.  
\end{DefNoNum} 

\begin{RemNoNum} 
The reference homeomorphism for the deformation is 
to fix an isotopy class of the model.  
In fact, only its isotopy class is significant.  
\end{RemNoNum} 

The topology on the set of equivalence classes 
of deformations of  $C$  turns out 
to be not quite complicated by the local rigidity.  
To see this, we review 
a few topological properties of the space of  
representations of the group  $\Pi = \pi_1(N)$  in 
$\PSL_2(\complex)$.  
The space of such representations carries 
a natural algebraic topology.  
There is a canonical projection to the set of  
conjugacy classes 
$\Hom(\Pi, \PSL_2(\complex)) / 
\PSL_2(\complex)$.   
A small neighborhood of a conjugacy class represented by  
a holonomy representation of a cone-manifold  
will be well behaved.  

An orientation preserving isometry  $\varphi$ 
of  $\hyp^3$  can be represented by a 
matrix  $\varPhi$  in  $\SL_2(\complex)$.  
The complex length of  $\varphi$  is a 
twice of an appropriate branch of the log of an
eigenvalue of  $\varPhi$.  
It measures how much  $\varphi$  
translates  $\hyp^3$  with twist along an invariant geodesic.  
Since there are two choices to be made, 
we adopt the following convention.  
To each oriented meridional element  $m_j $,  we choose 
a complex length of  $\rho_C(m_j)$  by  
a multiple of its cone angle with the complex unit $\sqrt{-1}$,  
which we denote by  ${\cal L}_{m_j}(\rho_C)$.   
It is equal to  $ \alpha^j \sqrt{-1}$.  
Then orient the invariant 
geodesic  $\ell_C^j$  of  $\rho_C(m_j)$  so that  
the rotational direction of  $m_j$  is counter clockwise.  
There are two ways to continuously extend 
this assignment of a complex length to 
that of transformation represented by  
$\rho(m_j)$  where  $\rho$  is close to  $\rho_C$.  
Since  $\rho$  is close to  $\rho_C$,  the invariant geodesic 
$\ell^j$  of  $\rho(m_j)$  is also close to  $\ell_C^j$  and 
inherits an orientation.  
We choose the sign of the real part   
according to whether  $\rho(m_j)$  translates  
the point on  $\ell^j$  in positive or negative direction.  

This convention extends to a continuous map, 
which we denote by  ${\cal L}_{m_j}$,  defined on the  
component 
$\Hom_C(\Pi, \PSL_2(\complex))/\PSL_2(\complex)$  of 
$\Hom (\Pi, \PSL_2(\complex))/\PSL_2(\complex)$   
containing the conjugacy class of  $\rho_C$.  
Arranging these maps, we obtain 
\begin{equation*}
{\cal L}_m : \Hom_C(\Pi, \PSL_2(\complex))/\PSL_2(\complex)
\to \complex^n, 
\end{equation*} 
where
${\cal L}_m (\rho) = 
({\cal L}_{m_1}(\rho), \cdots, {\cal L}_{m_n}(\rho))$.  
\bigskip

\begin{ThNum}[Hodgson-Kerckhoff's 
Local Rigidity \cite{HodgsonKerckhoff}]\label{LocalRigidity}
If the cone angles assigned to the components of  
$\varSigma$  in  $C$  all are positive and  $\leq 2 \pi$,  
then  ${\cal L}_m$  is a local diffeomorphism near the 
conjugacy class represented by   $\rho_C$.  
\end{ThNum}

In particular, the complex dimension of  
$\Hom(\Pi, \PSL_2(\complex)) / 
\PSL_2(\complex)$ 
near the class represented by  $\rho_C$  is 
equal to  $n$,  
the number of components of  $\varSigma$, and  
$\rho_C$  represents a smooth point.  

Removing a small tubular neighborhood of   $\varSigma$  from   $C$,  
we obtain a compact hyperbolic manifold with boundary.  
It's nearby deformations supported on all 
but a small neighborhood of the boundary are parameterized 
by holonomy representations  \cite{ThurstonNote}.  
The parameter in terms of complex lengths is 
more geometrically described in \cite{HodgsonKerckhoff, 
ThurstonNote}.  
For each value near  ${\cal L}_m(\rho_C)$, 
there is a unique way to fill the boundary  
by the hyperbolic Dehn filling theory 
\cite{ThurstonNote}. 
If the value of  ${\cal L}_{m}$  is 
in purely imaginal part    
$\Theta = (\Im \complex)^n \subset \complex^n$, 
then the filled resultant is a hyperbolic cone-manifold
with the same topological type but with perturbed 
cone angles. 
Hence the local rigidity together with 
the hyperbolic Dehn surgery filling imply 
the unique existence of a small deformation of  $C$  if 
the perturbation of cone angles is small enough.  
Namely the possible range of perturbation is open.  

When we start with a noncompact 3-cone-manifold 
without singularity, such as  $C_{comp}$  
then the existence of a small deformation is nothing but the 
conclusion of the hyperbolic Dehn filling theory 
for complete manifolds with a specified direction.  
In particular, there is a unique 
small deformation with perturbed cone angles 
for this case also.  
One can summarize the conclusion by  

\begin{CorNum}\label{SmallDeformationExist} 
Let  $C$  be an orientable hyperbolic 
3-cone-manifold of finite volume 
so that the cone angles assigned to 
the components of  $\varSigma$  all are at most  $2\pi$,   
possibly some or all of them are zero.  
Then there is a unique small deformation of  $C$  with 
perturbed cone angles if the perturbed cone angles 
are close enough to the initial ones. 
\end{CorNum}

\begin{RemNoNum}  
A continuous path  
on the space of representations whose image by  
${\cal L}_m$  is contained in purely imaginal part  
does not always correspond to  
the deformations of a cone-manifold in the full range.  
\end{RemNoNum}


\subsection{Pointed Hausdorff-Gromov Topology} 

We review the pointed geometric convergence 
of complete metric spaces due to Gromov 
in \cite{Gromov},  which generalizes the idea of  
Hausdorff topology on the set of all compact subset in 
a complete metric space.  

\begin{DefNoNum}
Let  $(X, x_0)$  and  $(Y, y_0)$  be pointed complete 
metric spaces.  
A relation  $\text{R} \subset X \times Y$  is an 
$\varepsilon$-approximation between  $(X, x_0)$  and  
$(Y, y_0)$  if  
\begin{enumerate} 
	\item
	there is  $y \in Y$  such that  $x_0 \text{R} y$  and  
		$d_Y(y_0, y) < \varepsilon$, 
	\item
	there is  $x \in X$  such that  $x \text{R} y_0$, 
		and $d_X(x_0, x) < \varepsilon$,   
	\item 
	$pr_X( \text{R} \cap (\ball_{1/{\varepsilon}}(X, x_0) 
		\times \ball_{1/{\varepsilon}}(Y, y_0)) 
		= \ball_{1/{\varepsilon}}(X, x_0)$,   
	\item 
	$pr_Y( \text{R} \cap (\ball_{1/{\varepsilon}}(X, x_0) 
		\times \ball_{1/{\varepsilon}}(Y, y_0)) 
		= \ball_{1/{\varepsilon}}(Y, y_0)$, and 
	\item 
	for any  $x, x' \in \ball_{1/{\varepsilon}}(X, x_0)$  and 
		$y, y' \in \ball_{1/{\varepsilon}}(Y, y_0)$  with 
		$x \text{R} y$  and  $x' \text{R} y'$, \par 
		we have 
		$\vert d_X(x,x') - d_Y(y,y') \vert < \varepsilon$. 
\end{enumerate}
\end{DefNoNum}

\begin{DefNoNum}  
A sequence  $\{ (X_i, x_i) \}$  of pointed complete 
metric spaces is said to converge geometrically to  $(Y, y)$  if 
for any  $\varepsilon > 0$, there exists   $i_0$  such that 
there is an $\varepsilon$-approximation between 
$ (X_i, x_i)$  and  $(Y, y)$  for all  $i > i_0$.  
\end{DefNoNum} 

An interesting class  ${\cal C}$  of 
pointed complete metric spaces with respect 
to this convergence is the one whose member  $(X, x)$ 
has the property;  
\begin{equation}\label{MetricBallCompact}  
	\text{The closure of} \quad \ball_R(X, x)  \quad 
	\text{is compact for all} \;\; R > 0. 
\end{equation} 
It is shown in  \cite{Gromov}  that 
the geometrically convergent sequence  $\{ (X_i, x_i) \}$  
contained in the class  ${\cal C}$  
geometrically converges to 
a unique complete metric space up to isometry.  
Hence the geometric convergence defines 
a Hausdorff topology on the set of isometry classes of 
${\cal C}$.  
We call it the pointed Hausdorff-Gromov topology.  

If furthermore we have some uniformity for 
the local structure of the metric spaces, 
then the set of isometry classes of such spaces
together with the pointed Hausdorff-Gromov 
topology becomes precompact.  
In particular, any sequence in the space contains 
a subsequence converging geometrically to a unique metric space.  
The following criterion by Gromov is useful.  

\begin{PropNum}[Gromov \cite{Gromov}]\label{GromovCriterion} 
Let  $\{ (X_i, x_i) \}$  be a sequence of pointed complete metric 
spaces with the property \eqref{MetricBallCompact} for all  $i$, 
then 
the sequence contains a geometrically 
convergent subsequence if and only if 
there is a subsequence  $\{ k \} \subset \{ i \}$  so that 
for any $R, \varepsilon > 0$,  
\begin{equation*} 
\min \# 
	\{ \varepsilon \text{\em {\rm -balls covering}} \,  
		\ball_R(X_k, x_k) \}
\end{equation*}
is bounded by some constant depending only on  $R$  and  
$\varepsilon$,   
where the minimum is taken over all 
$\varepsilon$-coverings of  $\ball_R(X_k, x_k)$.  
\end{PropNum} 

The class of metric spaces 
we are concerned with is quite restricted, but 
not only hyperbolic cone-manifolds.  
We will work also on euclidean cone-manifolds 
as their rescaling limits.  
 
\begin{NotNoNum} 
Fix a number  $L \leq -1$,  which will bound curvature.   
${\cal C}^{\theta}_{[L,0]}$  will be 
the set of pointed compact orientable 
cone-manifolds of constant sectional curvature $= K$,   
where 
$K \in [L, 0]$,   so that the cone angles assigned 
to the components of the singularity 
all are at most  $\theta$.  
${\cal C}^{\theta}_K$  is a subset of 
${\cal C}^{\theta}_{[L,0]}$  consisting of 
cone-manifolds with a particular curvature constant $K$.  
\end{NotNoNum} 

\begin{RemNoNum} 
Each member of  
${\cal C}^{\theta}_{[L,0]}$  is compact and hence has the 
property  \eqref{MetricBallCompact}.  
\end{RemNoNum}

Choose   $C \in {\cal C}^{2 \pi}_{K}$  and  
$x \in C$  not lying on  $\varSigma$.  
Here we bound cone angles by  $2\pi$. 
The set of points in  $C$  which admit at least two 
shortest geodesic paths to  $x$  in  $C$  is called 
a cut locus of  $C$  with respect to the 
reference point  $x$.  
The cut locus is a connected geodesic cell complex 
in  $C$.  
The complement of the cut locus consists of shortest 
rays to  $x$,  
\begin{equation*} 
P_x = \{ y \in C \, \vert \, y \; \text{admits the unique shortest 
path to} \; x \}, 
\end{equation*} 
and is called a 
Dirichlet fundamental domain of  $C$  about  $x$.  

\begin{LemNum}
The Dirichlet fundamental domain  $P_x$  of 
$C \in {\cal C}^{2\pi}_K$  about  $x$  is isometrically 
realized as an interior of 
a starshaped polyhedron in 
the simply connected 3-dimensional space  ${\bold H}_K$  of 
constant curvature $=K$.  
The closure is a starshaped polyhedron.  
Furthermore if  $C \in {\cal C}^{\pi}_K$,  
that is to say, cone angles all are $\leq \pi$,  then 
$P_{x}$  is convex.  
\end{LemNum} 

\begin{pf} 
Obvious. 
\end{pf}

We call this embedded compactified polyhedron  
a Dirichlet polyhedron of  $C$   about  $x$,  
and denote it again by  $P_x$. 
Namely,  $P_x$  stands for an open dense subset  
in  $C$,   and simultaneously, 
a compact polyhedron in  ${\bold H}_K$.  

\begin{PropNum}\label{GeometricLimitExist} 
Fix a curvature bound  $L \leq -1$.  
Any sequence  $\{ C_i \} \subset  {\cal C}^{2 \pi}_{[L, 0]}$  
with reference points  $x_i \in C_i$  
contains a subsequence converging geometrically to 
a complete metric space  $(C_*, x_*)$.  
\end{PropNum}

\begin{pf} 
Let  $(C, x_0)$  be a pointed compact orientable 
3-cone-manifold with constant sectional 
curvature $= K$  where  $L \leq K \leq 0$.  
If  $x_0 \not\in \varSigma$,  then 
\begin{align*}
\min \# 
	\{ \varepsilon \text{\rm -balls covering } \ball_R(C, x_0) \} 
		\leq &  
\min \# 
	\{ \varepsilon \text{\rm -balls covering } 
			\ball_R(P_{x_0}, x_0) \} \\
		\leq & 
\min \# 
	\{ \varepsilon \text{\rm -balls covering } 
		\ball_R({\bold H}_{K}, x_0) \}  \\ 
		\leq & 
\min \# 
	\{ \varepsilon \text{\rm -balls covering } 
		\ball_R({\bold H}_{L}, x_0) \},  
\end{align*}
where the minimum is taken all over  $\varepsilon$-ball 
coverings of the target.  
The last bound depends only on  $R$  and  $\varepsilon$,   
and not on any particular  $C$  or  $K$.  
If  $x_0 \in \varSigma$,  the same bound actually works 
by choosing the center of a Dirichlet polyhedron 
near  $x_0$  but not on $\varSigma$. 
Thus the result follows from 
Gromov's Criterion.  
\end{pf}

\begin{DefNoNum}
We call  $(C_*, x_*)$   a geometric 
limit of  $\{ (C_i, x_i) \}$.  
\end{DefNoNum} 

\begin{RemNoNum} 
If  $C_*$  is compact, the isometry class of  $C_*$  
does not depend on the choice of the reference points 
$x_i \in C_i$  in the sequence.  
\end{RemNoNum}

\begin{RemNoNum} 
The set of equivalence classes of deformations of 
an orientable hyperbolic 3-cone-manifold  $C$  
with cone angle  $\leq 2\pi$  carries a topology 
well described by the local rigidity.  
Assigning the isometry class of a deformation to 
each deformation, we get a map to the set of isometry classes 
of metric spaces in  ${\cal C}^{2\pi}_{-1}$  together 
with the pointed Hausdorff-Gromov topology.  
It is not quite hard to see that 
this map is continuous.  
\end{RemNoNum}


\section{Cone Angles $\leq \pi$}

In this section, we discuss three relative 
constants for hyperbolic 3-cone-manifolds with 
cone angles  $\leq \pi$,  
two of which dominate the local structure 
away from the singularity, and 
the other one of which is related with 
the geometry of a tubular neighborhood of the singularity.  
The angle condition  ``$\leq \pi$'' does  
not explicitly appear, but instead the fact 
that the Dirichlet polyhedron is convex, which is 
a conclusion of the angle condition, 
will be used often. 
We also discuss how cusp opening deformations occur 
locally.


\subsection{Thin Parts}

The constant in the following lemma is to claim 
that the injectivity 
radius decreases uniformly, like in the 
hyperbolic manifolds, away from the singularity.  

\begin{LemNum}\label{ThinPartFar} 
Fix a curvature bound  $L \leq -1$.  
Given positive numbers  $D, \; I, \; R > 0$,   
there is a constant  $U(D, I, R, L) > 0$  such that 
if  $C \in {\cal C}^{\pi}_{[L, 0]}$,  
$x \in C$  with  $d(x, \varSigma) \geq D$  
and  $\inj_x C \geq I$,  
then  
\begin{equation*} 
	\inj_y C \geq U(D, I, R, L)  
\end{equation*} 
for any  $y \in  C$      
with  $d(y, \varSigma) \geq D$  and  $d(y, x) \leq R$.  
\end{LemNum}

\begin{pf} 
Assume that there are no such uniform bounds.  
Then there is a sequence of cone-manifolds  
$\{ C_i  \} \subset C^{\pi}_{[L, 0]}$  and 
points  $ x_i, y_i \in C_i $  
such that  
\begin{enumerate} 
	\item $d(x_i, \varSigma_i), \; d(y_i, \varSigma_i) \geq D$, 
	\item $\inj_{x_i} C_i  \geq I$, 
	\item $d(y_i, x_i) \leq R$  \; and 
	\item $\inj_{y_i} C_i \leq 1/i$.  
\end{enumerate} 
Take a Dirichlet polyhedron  $P_{y_i}$  of  $C_i$  
about  $y_i$  in  ${\bold H}_{K_i}$,  
where  $K_i$  is a curvature constant of  
$C_i$.  
There are points  $p_i, q_i$  on  $\partial P_{y_i}$  which 
are identified in  $C_i$  and 
attain the shortest distance to  $y_i$  from  $\partial P_{y_i}$.  
The union of these paths forms a homotopically 
nontrivial shortest loop  $\ell_i$  in  $C_i$  
based at  $y_i$.  

$p_i$  and  $q_i$  are on the interior of the faces of  $P_{y_i}$  
respectively.  
Since the cone angles are  $\leq \pi$ and hence  
the Dirichlet polyhedron  $P_{y_i}$ is convex,  
$P_{y_i}$  is bounded by the extension of two faces 
which support  $p_i$  and  $q_i$.  
If the faces tend to be parallel as  $i \to \infty$,  
then the volume of  $\ball_{R+I}(C_i, y_i)$  
approaches zero.  
This is a contradiction since it must contain 
the ball  $\ball_{I}(C_i, x_i)$  whose volume 
admits nonzero lower bound by (2).  

If not,  
$\ell_i$  meets at  $y_i$  with 
angle uniformly away from  $\pi$.  
Let us lift  $\ell_i$  to a geodesic 
segment  $s_i$  in  $\hyp_{K_i}$  based at  $y_i$  such 
that  $p_i$  is the middle point.  
Then  $\rho_i(\ell_i)$  acts  ${\bold H}_{K_i}$   
by either a loxodromic 
(translation with twist) motion or an elliptic rotation. 
In both cases,  
the orbit of  $s_i$  by 
the action of a group generated by  $\rho_i(\ell_i)$  
forms a piecewise geodesic immersed line rounding 
around the axis of  $\rho_i(\ell_i)$.  
Since the corner of this line at the orbit of  $y_i$  
has an angle uniformly away from  $\pi$  with respect to  $i$, 
and since length of  $s_i$,  which equals the 
length of  $\ell_i$,  approaches zero when  $i \to \infty$, 
the immersed line must squeeze onto 
or into the axis of  $\rho_i(\ell_i)$  according to 
whether  $\rho_i(\ell_i)$  is loxodromic or elliptic.  
In particular in both cases, the axis of  $\rho_i(\ell_i)$  
becomes close to  $y_i$  when  $i \to \infty$.  

If there is a subsequence  
$\{ k \} \subset \{ i \}$  such that  
$\rho_k(\ell_k)$  all are loxodromic, 
then the translation distance becomes also short.  
The path joining  $p_k$  and  $q_k$  in  ${\bold H}_{K_k}$  
is equivariant tiny homotopic to the unit translation segment on 
the axis of  $\rho_k(\ell_k)$.  
This equivariant homotopy must induce a tiny homotopy in  $C_k$  
because of (1).  
Hence we obtain a very short closed geodesic in  $C_k$  
near  $y_k$.  
Then if we choose a new reference point  $z_k$ on 
this geodesic, the Dirichlet polyhedron 
about  $z_k$  will be bounded by almost parallel 
faces.  
This is a contradiction as before.  

In the other case,  $\rho_k(\ell_k)$  all but 
finitely many exceptions are elliptic.  
Then the axis of 
an elliptic element comes close to  $y_k$.   
Thus the length shortening tiny homotopy of  $\ell_k$  
must hit the singularity by 
Lemma~\ref{LengthShorteningHomotopyExist}, 
though the hit singularity may not be the axis 
of  $\rho_k(\ell_k)$. 
This contradicts (1).  
\end{pf}


\subsection{Local Margulis} 

The Margulis lemma for hyperbolic manifolds 
states that there is a universal 
constant depending only on the dimension 
which dominates the geometry and topology of thin parts.  
The cone-manifold admits in fact no such universal constant, 
however 
we may expect its relative version away from the 
singularity.  
The next lemma establishes that there is a Margulis like 
constant to control the geometry and topology of not absolute but 
relatively thin part with respect to the injectivity radius. 
We call it a local Margulis constant.  

\begin{LemNum}\label{LocalMargulis} 
Given positive numbers  $D, \, R > 0$,  
there is a constant   $V(D, R)$  such that 
if  $C \in {\cal C}^{\pi}_{-1}$,  
$d(x, \varSigma) \geq D$  and  $\inj_x C \leq V(D, R)$,  
then 
$(\ball_{R\cdot \inj_x C}(C, x), x)$  is homeomorphic to   
$(\ball_{R}(E, e), e)$  for some 
noncompact euclidean manifold  $E$  with 
$\inj_e E = 1$.  
\end{LemNum}

\begin{pf}
Assume that the conclusion is not true.  
Then there is a sequence of cone-manifolds  
$\{ C_i \} \subset  {\cal C}^{\pi}_{-1}$  and 
points  $x_i \in  C_i $  such that  
\begin{enumerate} 
	\item $d(x_i, \varSigma_i) \geq D$  and 
	\item $\inj_{x_i} C_i \leq 1/i$, but 
	\item $(\ball_{R \cdot \inj_{x_i}C_i} (C_i, x_i))$  never be   
		homeomorphic to  
		$(\ball_{R}(E, y), y)$  
		for some euclidean manifold with  $\inj_y E = 1$.  
\end{enumerate} 
Then 
\begin{equation*} 
	R \leq i D \leq D / \inj_{x_i}C_i 
\end{equation*} 
for  $i$  large enough, and we have 
\begin{equation*}
	R \cdot \inj_{x_i}C_i \leq D 
\end{equation*} 
for sufficiently large  $i$.   
Hence 
$(\ball_{R \cdot \inj_{x_i}C_i}(C_i, x_i), x_i)$  
is a subset of   
$(\ball_{D}(C_i, x_i), x_i)$  for 
sufficiently large  $i$.  
Notice that  $\ball_{D}(C_i, x_i)$  is nonsingular, 
so is  $\ball_{R \cdot \inj_{x_i} C_i} (C_i, x_i)$.  

Multiplying  $1/{\inj_{x_i}C_i}$  on the metric of  $C_i$,  
we obtain a cone-manifold  $\overline{C}_i$  
of constant curvature 
$ = - (\inj_{x_i}C_i)^2 \geq -1$  such that  
$\inj_{\overline{x}_i}\overline{C}_i = 1$.  
Hence we have a sequence of compact orientable 
cone-manifolds  $\{ \overline{C}_i \}$  
in  ${\cal C}^{\pi}_{[-1, 0]}$.  
Then by Proposition~\ref{GeometricLimitExist}, 
there is a subsequence  $\{ k \} \subset \{ i \}$  
so that  $\{ (\overline{C}_k, \overline{x}_k) \}$  
converges geometrically to a complete metric space  
$(\overline{C}_*, \overline{x}_*)$. 

The limit  $\overline{x}_*$  of reference points  
$\{ \overline{x}_k \}$  admits 
a neighborhood which is a limit of balls of radius $1$  
whose curvature tend to zero.  
On the other hand, an euclidean ball of 
radius  $1$  could be a geometric limit of this sequence.  
Hence by the uniqueness of the geometric limit, 
$\overline{x}_*$  admits an euclidean ball neighborhood.  
This point will be a reference point for the other part. 

To see a neighbor structure of the other part of 
$\overline{C}_*$, 
fix a constant  $R_1 > 0$  and 
choose any   $\overline{y}_* \in  \overline{C}_*$  with  
$d(\overline{y}_*, \overline{x}_*) \leq R_1$,  
then it is a limit of points  
$\{ \overline{y}_k \in \overline{C}_k \}$  
with, say,  $d(\overline{y}_k, \overline{x}_k) \leq 2R_1$.  
By rescaling, we have 
\begin{equation*} 
d(\overline{x}_k, \overline{\varSigma}_k) 
	\geq D/\inj_{x_k}C_k \geq kD, 
\end{equation*} 
and moreover, 
\begin{equation*} 
d(\overline{y}_k, \overline{\varSigma}_k) 
	\geq kD-2R_1 \geq D  
\end{equation*} 
for sufficiently large  $k$.  
Then by Lemma~\ref{ThinPartFar}, 
$\overline{y}_k$  admits a ball neighborhood of 
radius  $\geq U(D, 1, 2R_1, -1)$.   
This radius bound does not depend on  $k$,  
and hence a point  $\overline{y}_*$  
admits a neighborhood which is a limit of balls 
of uniformly bounded radius whose curvature tend to zero.   
Thus  $\overline{x}_*$  admits 
an euclidean ball neighborhood. 
Now, since  $R_1 > 0$  was arbitrary, a point with 
an arbitrary long distance from  $\overline{x}_*$  admits an 
euclidean ball neighborhood.  
This shows that 
$\overline{C}_*$  is an euclidean manifold 
without singularity.  
Moreover  $\inj_{\overline{x}_*} \overline{C}_* = 1$  and 
$\overline{C}_*$  is certainly noncompact.    

Letting  $(E, e) = (\overline{C}_*, \overline{x}_*)$, 
we will see that  $E$  has the property in the claim.  
Triangulate  $E$  by geodesic tetrahedra 
of uniform size and shape at least in 
a large compact set.   
Since it is a geometric limit, 
we may choose an approximate map  $\varphi_k$  from 
the $0$-skeleton of 
a large compact set of  $E$  containing 
$\ball_{D/ \inj_{x_k} C_k} (E, e)$  to  
a large compact set of  $\overline{C}_k$  containing   
$\ball_{D/ \inj_{x_k}C_k}(\overline{C}_k, \overline{x}_k)$  
for sufficiently large  $k$,  where 
4 vertices spanning an oriented simplex in  $E$  are mapped 
to 4-vertices spanning a simplex with the same 
orientation in the image.  
Then  $\varphi_k$  admits an obvious 
piecewise linear extension, 
which we denote again by  $\varphi_k$,  over 
a large compact set of $E$.  
$\varphi_k$  might be locally branched along edges or vertices, 
however it will be a homeomorphism for further sufficiently 
large  $k$,  because otherwise, 
$\varphi_k$'s  would not be accurate approximations. 
In particular, the restriction of  $\varphi_k$  
for sufficiently large  $k$  induces a homeomorphism of  
$(\ball_{R}(E, e), e)$  to 
$(\ball_{R}(\overline{C}_k, \overline{x}_k), \overline{x}_k)$  
after some tiny smoothing  
and hence to  
$(\ball_{R \cdot \inj_{x_k}C_k}(C_k, x_k), x_k)$.  
This is a contradiction.  
\end{pf}

\begin{RemNoNum} 
A homeomorphism 
between  $(\ball_{R \cdot \inj_x C}(C, x), x)$  and  
$(\ball_R(E, e), e)$  can be chosen by the composition 
of an approximation and a rescaling, 
which is an almost equi-expansive map centered at  $x$.  
\end{RemNoNum}

To see more topological structures of  $\ball_R(E, e)$,  
recall that a noncompact euclidean manifold is 
a quotient of the euclidean space  $\euc^3$  by 
a lattice  $\Gamma$  in  $\Isom_+ \euc^3$.  
$\Gamma$  is isomorphic to either  $\{0\}$,  
$\integer$, $\integer \oplus \integer$  or the 
fundamental group of the Klein bottle. 
In particular, any two elements in  $\Gamma$  are 
either mutually commute or anti-commute.  

\begin{LemNum}\label{VirtuallyAbelian} 
Let  $E$  be a noncompact euclidean manifold 
with a reference point  $e \in E$  such that  $\inj_e E = 1$.  
Then for any  $R > 0$,  the image of 
$\pi_1(\ball_{R}(E, e),e) \to \pi_1(\ball_{2R}(E, e),e)$  
induced by the inclusion is 
virtually abelian.  
\end{LemNum} 

\begin{pf} 
$\pi_1(\ball_R(E,e),e)$  is generated by 
geodesic loops of length  $\leq 2R$  which are  
not smooth only at the reference point  $e$.  
Choose any two such loops  $\ell_1, \ell_2$.  
Since they are commute or anti-commute in  $\pi_1(E)$, 
the lift of a loop represented by    
$\ell_1 \ell_2 \ell_1^{-1} \ell_2^{\varepsilon}$,  
where  $\varepsilon = 1 \, \text{or} \, -1$,   
encloses a rectangular knot in  $\euc^3$  where 
the edge length is at most  $2R$.  

The vertices are four lifts of the base point  $e$.  
The preimage of  
$\ball_{2R}(E, e)$,  which has twice in radius, 
contains a union of four balls of radius  $2R$  whose center 
lie on the vertices of the knot.  
Then the knot bounds a disk in the union of these four balls  
and hence $\ell_1 \ell_2 \ell_1^{-1} \ell_2^{\varepsilon}$  
becomes homotopic to zero in   $\ball_{2R}(E, e)$.  
This shows that the image is generated by finitely many  
$\ell$'s 
where the generators are mutually commute or 
anti-commute.  
Such a group is virtually abelian.  
\end{pf}


\subsection{Geometry of Tubes}

An abstract model for an 
equidistant tubular neighborhood 
of a singular geodesic in a cone-manifold  
will be useful for estimating several quantities.  
We call it a tube and discuss its geometry. 

\begin{NotNoNum} 
$T_{\sigma, \delta, \theta, \tau}$  will be 
an equidistant 
tubular neighborhood of a singular component 
with radius $\delta$  in 
a hyperbolic 3-cone-manifold where the length of 
a singular axis  
$= \sigma$,  the cone angle $= \theta$  and the 
twisting factor  $= \tau$.   
\end{NotNoNum}

These four parameters determine the isometry 
class of  $T_{\sigma, \delta, \theta, \tau}$.  
The boundary  
$\partial T_{\sigma, \delta, \theta, \tau}$  
carries an induced euclidean structure.  
A canonical rectangular 
fundamental domain of 
$\partial T_{\sigma, \delta, \theta, \tau}$  by 
the meridional direction and its vertical direction 
has magnitude   
\begin{equation}\label{Rectangle}
\theta \sinh \delta \times \sigma \cosh \delta.  
\end{equation} 
The surface area and volume of a tube 
depend only on the first three parameters,   
\begin{align} 
\label{TubeArea}
\area \partial T_{\sigma, \delta, \theta, \tau} & = 
	\theta \sigma \sinh \delta \cosh \delta, \\
\label{TubeVolume}
\volume T_{\sigma, \delta, \theta, \tau} & = 
	\frac{1}{2} \, \theta \sigma \sinh^2 \delta.   
\end{align} 
In particular,  
\begin{equation}\label{AreaVsVolume}
\frac{\area \partial T_{\sigma, \delta, \theta, \tau}} 
	{\volume T_{\sigma, \delta, \theta, \tau}}
	=  2 \coth \delta.   
\end{equation} 
The twisting factor  $\tau$  does not 
affect the surface area and volume in fact, 
but does affect the euclidean structure 
of the boundary.  
Since the gluing of the right and left edges has 
no twisting factor, we have by \eqref{Rectangle}  
\begin{equation}\label{InjBound}
\inj \partial T_{\sigma, \delta, \theta, \tau} 
	\leq  \theta \sinh \delta. 
\end{equation} 
The second factor 
``$\sigma \cosh \delta$''  of  \eqref{Rectangle}  
does not say much about injectivity radii.   
\medskip

We will derive two implications from these quantities. 
They are about rank 2 cusp opening  
deformations of a tube, 
and the comparison of intrinsic and extrinsic 
injectivity radii of a point on the boundary of 
a tube embedded in a cone-manifold.  

First of all, 
consider a sequence  $\{ T_i \}$  of tubes which converges 
geometrically to a rank 2 cusp neighborhood by 
taking reference points on the boundaries. 
There are 
essentially two different ways.  
The simplest one in one way can be seen by 
setting 
\begin{align*} 
\theta_i & = 1/\sinh \delta_i,  \\
\sigma_i & = 1/\cosh \delta_i 
\end{align*} 
and making  $\theta_i \to 0$  when  $i \to \infty$.  
Then  $\delta_i \to \infty$  and  $\sigma_i \to 0$.  
There are an elliptic and a loxodromic 
elements in  $\PSL_2(\complex)$  which generate  $T_i$  for 
each  $i$.  
In this deformations, 
they both 
approach parabolic 
elements which generate a rank 2 cusp when  $i \to \infty$.  

The other way involving a twisting factor  $\tau$  was 
discussed in  \cite{ThurstonNote}.  
To see one simple example, let  $\theta$  be a fixed 
positive constant and set  
\begin{equation*} 
\sigma_i \sinh \delta_i \cosh \delta_i = 1 /\theta.  
\end{equation*} 
Choose  $\tau_i$  so that  $\{ \partial T_i \}$  
stays in a compact set in the moduli space of 
euclidean tori, we obtain a cusp opening family.  
In this case,  $\delta_i \to \infty$  and  $\sigma_i \to 0$  also.  
However, the elliptic elements associated to  
$T_i$'s  diverge in  $\PSL_2(\complex)$, and 
a group generated by a loxodromic element 
converges geometrically (but not algebraically)  
to a rank 2 parabolic subgroup generating 
a cusp.  

\begin{LemNum}\label{CuspOpenning}
Let  $\{ T_i \}$  be a sequence of tubes which 
converges geometrically to a rank 2 cusp neighborhood. 
Suppose that  $\area \partial T_i$  is constant and 
that  $\{ \partial T_i \}$  stays in a compact set 
in the moduli space of euclidean tori.  
If  $\theta_i > 0$  is bounded away from zero, 
then the lengths of curves on  $\partial T_i$'s 
bounding a disk in  $T_i$  diverge as  $i \to \infty$.  
\end{LemNum} 

\begin{pf} 
Since the boundaries have bounded 
geometry and  $T_i$  approaches a rank 2 cusp 
neighborhood, 
$\delta_i \to \infty$.  
Then since the length of a curve bounding 
a disk in  $T_i$  is  $\theta_i \sinh \delta_i$, 
it must diverge if  $\theta_i$  is bounded 
away from zero. 
\end{pf} 

Secondly, consider a tube   
$T^K_{\sigma, \delta, \theta, \tau}$  
with constant sectional curvature  $= K$,  
where  $K$  lies in  $[L, 0]$  and is not 
necessarily  $-1$.  
The modified  formulas of the above ones for  
$T^K_{\sigma, \delta, \theta, \tau}$  
can be established 
according to the value of  $K$,   however 
its limiting behaviors caused by 
the limiting behavior of  $\sigma, \delta, \theta$  
are essentially the same.  

Regarding a tube as a model of an equidistant tubular 
neighborhood  ${\cal T} \subset C$  of some 
component of  $\varSigma$,  
we compare the injectivity radius of $\partial {\cal T}$  
and injectivity radius of  $C$  at a point 
on  $\partial {\cal T}$.  
Notice that since  $\partial {\cal T}$  is an 
euclidean torus,  $\inj_x \partial {\cal T}$  does not 
depend on the choice of  $x \in \partial {\cal T}$  and 
is equal to  $\inj \partial {\cal T}$.  

\begin{LemNum}\label{InjectivityComparison} 
Fix a curvature bound  $L \leq -1$.  
Given positive numbers  $D, I > 0$,  there is 
a constant  $W(D, I, L)$  such that 
if  $C \in {\cal C}^{\pi}_{[L, 0]}$,  
${\cal T} \subset C$  is an equidistant tubular neighborhood 
of a component of  $\varSigma$    
with  $\radius {\cal T} \geq D$  and  
$\inj \partial {\cal T} \leq I$,  then 
\begin{equation*} 
	(\inj_wC \leq ) \; 
	\inj \partial {\cal T} \leq W(D,I,L) \cdot \inj_w C 
\end{equation*} 
for any  $w \in \partial {\cal T}$, where 
the inequality in $( \quad )$ on the left hand side is obvious.  
\end{LemNum}

\begin{pf} 
Let  $K'$  be a curvature constant of  $C$.  
Multiplying  $1/\inj \partial {\cal T}$  on 
the metric of  $C$,    
we obtain  $\overline{C}$  of 
constant curvature 
$K = K'(\inj \partial {\cal T})^2 \geq I^2L$.   
Then  ${\cal T}$  
becomes a tubular neighborhood  $\overline{{\cal T}}$  with 
$\radius \overline{{\cal T}} \geq D/I$ and also 
$\inj \partial\overline{{\cal T}} = 1$.  

Let 
$T^K_{\sigma, \delta, \theta, \tau}$  be 
a tube isometric to  $\overline{{\cal T}}$, 
where  $\delta = \radius \overline{{\cal T}} \geq D/I$, 
and 
consider a natural inclusion of tubes 
\begin{equation*} 
T^K_{\sigma, \delta-D/2I, \theta, \tau} \subset 
T^K_{\sigma, \delta, \theta, \tau}.  
\end{equation*} 
Notice that the radius of the included tube is  
$\geq D/2I$.  
$T^K_{\sigma, \delta, \theta, \tau}$  is 
a Riemannian manifold with boundary, and    
the quantity 
$\inj_x T^K_{\sigma, \delta, \theta, \tau}$  
for 
$x \in \partial T^K_{\sigma, \delta-D/2I, \theta, \tau}$   
can be considered as a positive periodic function 
in terms of  $\tau$,  
so that we have the minimum for fixed 
$\sigma, \delta, \theta, K$.  
\begin{equation*}
J^K_{\sigma, \delta, \theta} (D/I) 
	= \min_{\tau} \inj_x T^K_{\sigma, \delta, \theta, \tau}.  
\end{equation*}
Then 
$J^K_{\sigma,\delta,\theta}(D/I)$  becomes  
a function of  
$\sigma > 0, \delta \geq D/I, \theta > 0$,  
and  $0 \geq K \geq I^2L$.  
Since it approaches zero only when 
$\inj \partial T^K_{\sigma, \delta, \theta, \tau} \to 0$, 
it attains the positive minimum  $J(D, I, L)$  
when the variables run so that 
$\inj \partial T^K_{\sigma, \delta, \theta, \tau} = 1$.  
In particular,  
$\inj_{\overline{x}} \overline{{\cal T}} \geq J(D, I, L)$  
for any  $\overline{x} \in \overline{{\cal T}}$  with  
$d(\overline{x}, \partial \overline{{\cal T}}) = D/2I$.  

Set a constant  
$U$  by  $U(D/2I, J(D,I,L), D/2I, I^2L)$  in 
Lemma~\ref{ThinPartFar},  and 
choose any  $\overline{w} \in \partial \overline{{\cal T}}$.  
Then there is the unique nearest point  
$\overline{x}$  to  $\overline{w}$  in  
$\overline{{\cal T}}$  so that  
$d(\overline{x}, \overline{w}) = d(\overline{x}, \partial \overline{{\cal T}}) 
= D/2I$.  
Since  $d(\overline{x}, \overline{\varSigma}) \geq D/2I$, 
$\inj_{\overline{x}} \overline{C} \geq 
\inj_{\overline{x}} \overline{{\cal T}} \geq J(D, I, L)$  and 
$d(\overline{w}, \overline{\varSigma}) \geq D/2I$, 
we have a bound 
$U \leq \inj_{\overline{w}} \overline{C}$  
by Lemma~\ref{ThinPartFar}.  
Multiplying  $\inj \partial {\cal T}$ on the metric, 
we obtain
\begin{equation*}
U \cdot \inj \partial {\cal T} \leq \inj_{w} C.   
\end{equation*} 
The proof is done by letting  $W = 1/U$. 
\end{pf}


\section{Choosing Geometrically Convergent Sequences}

In this section, 
we choose a geometrically convergent sequence 
of deformations of a 
compact orientable hyperbolic 3-cone-manifold  $C$  
with singularity  $\varSigma$.  
This is to see what happens in the limit in the next 
two sections.


\subsection{Maximal Tubes}

We will often split the discussions in the later section 
according to whether 
the reference point is close to the singularity or not. 
For this, it is convenient to 
introduce the maximal tubular neighborhood of  
the singularity.  

\begin{DefNoNum}
The maximal tube  ${\cal T}$  about  $\varSigma \subset C$  is a 
union of open tubular neighborhoods  ${\cal T}^j$'s  such that 
\begin{enumerate} 
\item  each component  ${\cal T}^j \subset {\cal T}$  is 
	an equidistant tubular 
	neighborhood to the  $j$th component 
	$\varSigma^j \subset \varSigma$ of the 
	singularity, 
\item  among ones having the property (1),  
	the set of radii arranged in order of magnitude  
	from the smallest one is maximal in lexicographical order.  
\end{enumerate} 
\end{DefNoNum}

By the second condition, the maximal tube  ${\cal T}$  
about  $\varSigma$   
is uniquely determined. 
In fact, ${\cal T}$ can be constructed as follows.  
Since the components of  $\varSigma$  are isolated, 
a $\delta$-tubular neighborhood of  $\varSigma$  is 
a union of disjoint tubular neighborhoods of  
$\varSigma^j$'s  
if  $\delta$  is very small.  
Thicken it gradually.  
Then some of components contact at 
a particular moment.  
Stop the growth of the components involved 
in contacting, 
and thicken the others furthermore.  
We will have the second contact moment.  
Do the same again.  
Continue this process up to the terminal moment, 
and we finally obtain the maximal tube  
${\cal T}$.  

Denote by  $\partial {\cal T}^j$  an abstract 
boundary of  ${\cal T}^j$.  
It is an euclidean torus. 
The actual boundary  $\partial {\cal T}$  
of  ${\cal T}$  in  $C$
is a union of isometrically embedded tori 
with a finite number of contact points.


\subsection{Geometrically Convergent Sequences}

Recall that  $A = ( \alpha^1, \cdots, \alpha^n)$  
is an angle set of a compact orientable hyperbolic 
3-cone-manifold  $C$,  where 
$\alpha^j \leq \pi$  for all  $1 \leq j \leq n$.  
Choose a sequence of deformations  $(C_i, \xi_i)$  of  $C$  
such that  
\begin{enumerate} 
\item $\alpha^j_i \leq \pi$  for any  $i$,  and  
\item $\lim_{i \to \infty} \alpha^j_i$  exists and equals  $\beta^j$.  
\end{enumerate} 
To see what happens when  $i \to \infty$,  
we follow Thurston's strategy in \cite{SOK} in the 
next two sections, that is, 
to analyze its possible geometric limit using 
hyperbolic geometry and 3-manifold topology.  

Let  $c_i$  be the first contact point on 
$\partial {\cal T}_i$,   
that is to say, the point which admits two 
shortest path to  $\varSigma_i$  from  
$\partial {\cal T}_i$.  
Then by Proposition~\ref{GeometricLimitExist}, 
there is a subsequence  $\{ k \} \subset \{ i \}$  
so that  $\{ (C_k, c_k) \}$  converges geometrically to a 
complete metric space  $(C_*, c_*)$.  
We arrange more.  
Let  $f_k$  be the finest point on  $\partial {\cal T}_k$, 
that is to say,  the point on  $\partial {\cal T}_k$  
which attains the minimum among 
$\{ \inj_x C_k \vert x \in \partial {\cal T}_k \}$.   
By choosing a further subsequence  $\{ C_k\}$  
with the same letter, we may assume that   
$c_k$  and  $f_k$  stay on the components  
of  $\partial {\cal T}_k$  with constant reference numbers 
respectively. 
Namely,  $c_k$  always lies in $c$th component  
$\partial {\cal T}_k^c$  for any  $k$,  and 
$f_k$  does in  $\partial {\cal T}_k^f$,  
where  $c$  and  $f$  here 
represent the reference numbers for the components.  

The properties of the sequence so chosen is 
summarized in  \ref{Setting}  below.  
In fact, we only use the properties listed in \ref{Setting} 
in the next two sections.


\subsection{Properties}\label{Setting}

{\it 
A sequence of compact orientable hyperbolic 3-cone-manifold  
$\{ C_k \}$  with the angle set  
$A_k = (\alpha^1_k, \cdots, \alpha^n_k)$  
has the properties:   
\begin{enumerate} 
\item\label{SettingOne} 
	each  $C_k$  is a deformation of  $C$  with a reference 
	homeomorphism  $\xi_k : C \to C_k$,  
\item\label{SettingTwo} 
	$\{ (C_k, c_k) \}$  converges geometrically 
	to  $(C_*, c_*)$  when  $k \to \infty$, 
\item\label{SettingThree} 
	$\alpha^j_k \leq \pi$  for all  $1 \leq j \leq n$  and 
	any  $k$,    
\item\label{SettingFour} 
	$\alpha^j_k \to \beta^j$  when  $k \to \infty$, 
\item\label{SettingFive} 
	the first contact point  $c_k$  lies on 
	a component  $\partial {\cal T}_k^c$  with 
	a constant reference number  $c$, 
\item\label{SettingSix} 
	the finest points  $f_k$  lies on a component  $\partial {\cal T}_k^f$  
	with a constant reference number  $f$  and  
\item\label{SettingSeven} 
	there is a constant  $V_{max}$  such that  
$\volume C_k \leq V_{max}$.   
\end{enumerate}
}

\begin{RemNoNum} 
The sequence above is assumed to have only a geometric limit and 
the algebraic convergence with respect 
to the identification by  $\xi_k$  is not guaranteed.  
For instance, let  $C$  be a hyperbolic surface with 
homotopically nontrivial 
automorphism  $\varphi : C \to C$  which cannot be realized by 
an isometry, and define  $\xi_k : C \to C_k = C$  by  
\begin{equation*} 
\xi_k = 
	\begin{cases} 
	id, & \text{if  $k$  odd}, \\
	\varphi, & \text{if  $k$  even}. 
	\end{cases} 
\end{equation*}
Then the sequence  $\{ C_k \}$  converges geometrically to  
$C$  but the a sequence of holonomy representations  $\{ \rho_k \}$  
does not 
converge algebraically.  
\end{RemNoNum} 

\begin{RemNoNum} 
The property \eqref{SettingSeven}  is not a direct consequence 
of  Proposition~\ref{Volume<Bound}, however the bound above 
can be obtained by modifying its proof only   
a little because of the property \eqref{SettingThree}. 
See also the 
remark after  Proposition~\ref{Volume<Bound}.  
\end{RemNoNum}


\section{Thick Tube}

In this section, we study what happens to a geometric limit  
$C_*$  of a sequence  $\{ C_k \}$  of the deformations 
of  $C$  in \ref{Setting} when maximal tubes of the 
singularity are uniformly thick.  


\subsection{Brief Outline}

The underlying assumption throughout this section is  

\begin{AssOwnNum}\label{ThickTube}
There is a constant  $D_1 > 0$  such that 
\begin{equation*}
D_1 \leq \radius {\cal T}_k^j 
\end{equation*}  
for any  $1 \leq j \leq n$  and any  $k$.  
\end{AssOwnNum} 

Under the Assumption~\ref{ThickTube} above, 
we prove the following propositions in this section.  

\begin{PropNum}\label{ThickProp1}
Under the Assumption~\ref{ThickTube},   
there is a constant  $I_1 > 0$  such that 
\begin{equation*} 
I_1 \leq \inj_{f_k}C_k. 
\end{equation*} 
for any  $k$.  
\end{PropNum} 

The conclusion of this proposition is 
equivalent to the Assumption~\ref{ThickTube}  
since $\inj_{f_k} C_k \leq \min_j \radius {\cal T}^j_k$.  
The proof involves analysis of the local structure 
of cone-manifolds away from the singularity.  

Using the conclusion of Proposition~\ref{ThickProp1}, 
we show 

\begin{PropNum}\label{ThickProp2}  
Under the conclusion of Proposition~\ref{ThickProp1},   
$C_*$ is 
a hyperbolic 3-cone-manifold of finite volume 
homeomorphic to  $C$,  where 
some components of  $\varSigma$  possibly disappear and 
create cusps.  
\end{PropNum}

\begin{PropNum}\label{ThickProp3} 
Under the conclusion of Proposition~\ref{ThickProp1},  
a sequence  $\{ \rho_k \}$  of holonomy representations   
of  $\{ C_k \}$  contains a subsequence converging  
algebraically to 
the holonomy representation  $\rho_*$  
of  $C_*$  with respect to the identification by  $\xi_k$.  
\end{PropNum} 

\begin{DefNoNum} 
A sequence  $\{ C_k \}$  of deformations is 
said to converge strongly if 
it converges geometrically to a cone-manifold  $C_*$  
homeomorphic to  $C$  
and a sequence  $\{\rho_k\}$  of 
their holonomy representations  
converges algebraically to  $\rho_*$  with 
respect to the identification by  $\xi_k$.  
\end{DefNoNum} 

This definition is compatible for the existing 
one for discrete groups.

\begin{CorNum}\label{StrongLimit} 
Under the Assumption~\ref{ThickTube}, 
$\{ C_k \}$  contains a subsequence which 
converges strongly to a 
hyperbolic 3-cone-manifold  $C_*$  
homeomorphic to  $C$.  
If  $\beta^j > 0$  for  all $1 \leq j \leq n$, then 
$C_*$  is compact.  
\end{CorNum}

\begin{pf} 
This is a direct consequence of three 
propositions above. 
Suppose that the Assumption~\ref{ThickTube} is the case.  
By Proposition~\ref{ThickProp1}, 
the injectivity radius of the first contact 
points of maximal tubes are uniformly bounded 
away from zero. 
Then by Proposition~\ref{ThickProp2}, \ref{ThickProp3}, 
$C_*$  is a strong limit of  $\{ C_k \}$  after 
taking a subsequence.  
In particular, the angle set of  $C_*$  is 
equal to  $B = (\beta^1, \cdots , \beta^n)$.  
If  $\beta^j > 0$  
for all  $1 \leq j \leq n$,  then  $C_*$  admits 
no ends and hence is compact.  
\end{pf}


\subsection{Boundary of Tubes}

Supposing that the Assumption~\ref{ThickTube} is the case  
throughout this subsection, we  
prove Proposition~\ref{ThickProp1}

\begin{LemNum}\label{AreaBound}
There is a constant  $I_2$  such that  
\begin{equation*} 
\area \partial {\cal T}_k^j \leq I_2, 
\end{equation*} 
for any  $1 \leq j \leq n$  and any  $k$.  
\end{LemNum} 

\begin{pf} 
By the comparison \eqref{AreaVsVolume}  of the 
volume and the surface area of a tube, and by the 
volume bound \eqref{SettingSeven} in \ref{Setting}, we have 
\begin{equation*} 
\area \partial {\cal T}^j_k  = 
	2 \coth \radius {\cal T}_k^j \, \volume {\cal T}_k^j 
	\leq  2 V_{max} \coth D_1.   
\end{equation*} 
Let  $I_2$  be the last term.  
\end{pf}

\begin{pf*}{\it Proof of Proposition~\ref{ThickProp1}}
Assume contrarily that  
$\inj_{f_k} C_k \to 0$, 
and we will get a contradiction.  

Choose the first contact point  $p_k$  on 
$\partial {\cal T}_k^{f}$  on which  $f_k$  lies.  
$p_k$  may not be the absolute first contact 
point  $c_k$  since we require that  
$p_k$  lies on the component  $\partial {\cal T}_k^f$  
which might be different from  $\partial {\cal T}_k^c$.    
$p_k$  is the point where 
$\partial {\cal T}_k^f$  either meets the other component of 
$\partial {\cal T}_k$  or contacts itself.  
Very locally,  $p_k$  appears as a contact point of 
two components, one  $\partial {\cal T}_k^f$  from the left 
hand side   
and the other  $\partial {\cal T}_k^{f'}$  from the right hand 
side.  
The reference numbers  $f$  and  $f'$  might be 
the same.  

One obvious inequality is 
\begin{equation*}
(\inj \partial {\cal T}^{f}_k)^2 \leq \area \partial {\cal T}_k^{f} 
\; ( \leq I_2),    
\end{equation*}
and the right hand side of 
which is bounded by  $I_2$  by the above lemma.  
Let  $W$  be a constant  $W(D_1, \sqrt{I_2}, -1)$  in 
Lemma~\ref{InjectivityComparison}.  
Then since  
$\radius {\cal T}_k^{f} \geq D_1$  and 
$\inj \partial {\cal T}^{f}_k \leq \sqrt{I_2}$,  
we have 
\begin{equation}\label{RatioBound1} 
\inj \partial {\cal T}_k^f \leq 
	W \cdot\inj_{f_k} C_k 
		\; (\leq W \cdot {\inj_{p_k} C_k})    
\end{equation} 
by Lemma~\ref{InjectivityComparison}.  
If we regard  $p_k$  as a point on  $\partial {\cal T}_k^{f'}$,  
then since  $\radius {\cal T}_k^{f'} \geq D_1$  and  
$\inj \partial {\cal T}_k^{f'} \leq \sqrt{I_2}$,  again   
we have 
\begin{equation}\label{RatioBound2} 
\inj \partial {\cal T}_k^{f'} \leq 
	W \cdot {\inj_{p_k} C_k}.  
\end{equation} 
Hence 
$\ball_{W \cdot \inj_{p_k} C_k}(C_k, p_k)$  contains 
a homotopically nontrivial loop  
$\ell_1$  based at  $p_k$  on the 
left  $\partial {\cal T}_k^f$,  
and also  $\ell_2$  based at  $p_k$  on 
the right  $\partial {\cal T}_k^{f'}$.  

Let  $V(D_1, 2W)$  be a local Margulis constant with 
respect to  $D_1$  and  $2W$  in 
Lemma~\ref{LocalMargulis}.  
Since we assumed that  $\inj_{f_k} C_k \to 0$, 
\eqref{RatioBound1} and  \eqref{RatioBound2}  imply that 
\begin{equation*} 
(\inj_{p_k} C_k \leq) 
	\; \inj \partial {\cal T}^f_k \leq 
		V(D_1, 2W). 
\end{equation*} 
for sufficiently large  $k$.  
Then 
$(\ball_{2W \cdot \inj_{p_k} C_k}(C_k, p_k), p_k)$  is 
homeomorphic by an almost equi-expansive map 
to  $(\ball_{2W}(E, e), e)$  for some 
noncompact euclidean manifold  $E$  with 
$\inj_e E = 1$  
by Lemma~\ref{LocalMargulis} and the remark after that.  
Furthermore  
by  Lemma~\ref{VirtuallyAbelian},  the homomorphism,  
\begin{equation*} 
\pi_1(\ball_W(E, e), e) 
\to \pi_1(\ball_{2W}(E, e), e),   
\end{equation*} 
induced by the inclusion
have a virtually abelian image.  
Hence so does 
\begin{equation*} 
\pi_1(\ball_{W \cdot \inj_{p_k} C_k}(C_k, p_k), p_k) 
\to \pi_1(\ball_{2W \cdot \inj_{p_k} C_k}(C_k, p_k), p_k),   
\end{equation*} 
because of the choice of 
a homeomorphism of  
$(\ball_{2W \cdot \inj_{p_k} C_k}(C_k, p_k), p_k)$  we made.  
Thus the nontrivial loops $\ell_1, \ell_2$  
representing elements of  
$\pi_1(\ball_{W \cdot \inj_{p_k} C_k}(C_k, p_k), p_k)$  
are virtually commutative, in particular, in  $\Pi$  
for sufficiently large  $k$.  

On the other hand, consider the developed image near  $p_k$.  
$p_k$  lifts to a contact point of lifts of  
$\partial {\cal T}_k^f$  from the left hand side and 
that of 
$\partial {\cal T}_k^{f'}$  from the right hand side.  
In particular,  
$\rho_k(\ell_1)$  leaves the left hand lift of  
$\partial {\cal T}_k^f$  invariant, on the other hand, 
$\rho_k(\ell_2)$  leaves the right hand lift of 
$\partial {\cal T}_k^{f'}$  invariant.  
Hence their action on  $\hyp^3 \cup \sphere^2_{\infty}$   
do not have common fixed point at all, 
and they are not commutative 
even virtually in  $\rho_k(\Pi)$.  
This is a contradiction.  
\end{pf*}


\subsection{Geometric Limits}

We assume 
the conclusion of Proposition~\ref{ThickProp1} 
that the injectivity radius of the points 
on  $\partial {\cal T}_k$  is uniformly bounded 
from below by  $I_1 > 0$, and prove 
Proposition~\ref{ThickProp2}.  

\begin{LemNum}\label{CoreLengthBound} 
Under the conclusion of Proposition~\ref{ThickProp1}, 
if there is a constant  $D^j > 0$  such that  
$\radius {\cal T}^j_k \leq D^j$,  
then 
$\beta^j > 0$  and there is a constant  $S^j > 0$  
such that  
\begin{equation*} 
S^j \leq \length \varSigma_k^j.  
\end{equation*} 
for any  $k$.  
\end{LemNum}

\begin{pf} 
If  $\beta^j = 0$, then the $j$th component 
$\alpha_k^j$  of the angle set approaches zero as 
$k \to \infty$.  
Thus since  $\radius {\cal T}_k^j \leq D^j$,   
$\inj \partial {\cal T}_k^j \to 0$  
by  \eqref{InjBound}.  
This contradicts the conclusion of 
Proposition~\ref{ThickProp1}  since
\begin{equation*} 
(I_1 \leq) \; \inj_{f_k} C_k \leq \inj \partial {\cal T}_k^j.  
\end{equation*} 

If  $\{ \length \varSigma_k^j \}$  
contains a subsequence converging to  $0$, 
then since 
\begin{align*} 
\area \partial {\cal T}_k^j 
	& = \theta_k^j \, 
	\length \varSigma_k^j \, \sinh \radius {\cal T}_k^j \,
	\cosh \radius {\cal T}_k^j \\
	& \leq \pi \, \sinh D^j \cosh D^j 
		\length \varSigma_k^j,
\end{align*} 
by \eqref{TubeArea},  
$(I_1^2 \leq (\inj \partial {\cal T}_k^j)^2 \leq) \; 
\area \partial {\cal T}_k^j$  
can be arbitrary close to  $0$.  
This is again a contradiction. 
\end{pf}

\begin{LemNum}\label{NoncompactLimit}
$C_*$  is a hyperbolic 3-cone-manifold of finite volume  
possibly with compact singularity.  
\end{LemNum} 

\begin{pf} 
It is sufficient to show that 
each point  $x_* \in C_*$  admits a hyperbolic 
ball neighborhood possibly with singularity 
along a geodesic segment.  
The argument is quite parallel to that in Lemma~\ref{LocalMargulis}.  

Since  $I_1 \leq \inj_{f_k} C_k \leq \inj_{c_k} C_k$,  
the limit  $c_* \in C_*$  of the reference points  
$\{ c_k \}$  admits 
a neighborhood which is a limit of 
hyperbolic balls of uniformly bounded radii by  $I_1$.   
Hence it admits a hyperbolic ball neighborhood.  

Fix a constant  $R > 0$  
and choose any  $x_* \in  C_*$  with 
$d(x_*, c_*) \leq R$.  
If it is a limit of 
points  $\{ x_k \in C_k - {\cal T}_k \}$,   
then since  $d(x_k, \varSigma) \geq D_1$    
and we may assume  $d(x_k, c_k) \leq 2R$,  
$x_k$  admits a hyperbolic ball neighborhood of 
radius  
$\geq U(D_1, I_1, 2R, -1)$  by Lemma~\ref{ThinPartFar},  
where the radius bound does not depend on  $k$.  
Hence again  $x_*$  admits a hyperbolic open ball 
neighborhood.  

If   $x_*$  is a limit of points  $\{ x_k \in {\cal T}_k \}$,  
we may assume that  $\{ x_k \}$  is contained in 
a  component  ${\cal T}^x_k$  with a constant 
reference number  $x$     
by taking further subsequence if necessary.  
When  $\radius {\cal T}_k^x \to \infty$, 
since  $\partial {\cal T}^x_k$  does not 
degenerate by the conclusion of Proposition~\ref{ThickProp1}, 
we may assume that 
$d(x_k, c_k) \leq 2R,  
d(x_k, \varSigma_k^x) \to \infty$,  and   
$d(x_k, \varSigma_k) \geq D_1$.  
Then  
$x_k$  admits a hyperbolic ball neighborhood of 
radius  $\geq U(D_1, I_1, 2R, -1)$   
by Lemma~\ref{ThinPartFar}  
where the radius bound does not depend on  $k$, 
and hence so does  $x_*$.  
When  $\radius {\cal T}_k^x$  is bounded not only 
from below by the Assumption~\ref{ThickTube}  but 
also from the above,  
$\length \varSigma_k^x$  has a uniform lower bound away from zero
by Lemma~\ref{CoreLengthBound}, 
and each point within  ${\cal T}^x_k$  has a 
possibly singular ball neighborhood of uniform radius 
where the singularity occurs only along a geodesic segment.  
Hence  $x_*$  admits a hyperbolic ball 
neighborhood possibly with cone singularity along 
a geodesic segment.  

Since  $R > 0$  was arbitrary,  
the above argument shows that 
every point on  $C_*$  admits a hyperbolic 
ball neighborhood possibly with a cone singularity.  
The singularity appears only in the limit of 
${\cal T}_k^x$  whose radius is bounded.  
There are only finitely many such components.  
Moreover the length of a core of such a component 
is bounded since we have in general by \eqref{InjBound}  
\begin{equation*} 
	(I_1 \leq) \, \inj \partial {\cal T}_k^x 
	\leq \theta_k^x \sinh \radius {\cal T}_k^x, 
\end{equation*} 
so that the formula in \eqref{TubeArea} implies the estimate, 
\begin{equation*} 
	\length \varSigma_k^x = 
	\frac{\area \partial {\cal T}_k^x}
	{\theta_k^x 
		\sinh \radius {\cal T}_k^x \cosh \radius {\cal T}_k^x} 
	\leq \frac{I_2}{I_1 \cosh \radius D_1}.  
\end{equation*} 
Hence the singular set is compact.  
\end{pf}

By Corollary~\ref{EndToral},  
$C_*$  has finitely many toral ends.  
Choose disjoint horotoral neighborhoods of the 
ends of  $C_*$  so that  
the minimum  $I_3$  of the injectivity radii of 
$C_*$  at points on their boundaries is  $\leq I_1$.    
We let  $C_*^{cut}$  
be a compact hyperbolic cone-manifold 
with toral boundary obtained 
from  $C_*$  by truncating such 
cusp neighborhoods.  
We thus have 
\begin{equation*} 
	I_3 = \min \{\inj_x C_* \vert x \in \partial C_*^{cut} \} \leq I_1.  
\end{equation*} 

\begin{LemNum}
There is an approximate homeomorphism  
$\varphi_k : C_*^{cut} \to C_k$  for sufficiently large  $k$. 
\end{LemNum} 

\begin{pf} 
We just repeat the last paragraph in the proof 
of Lemma~\ref{LocalMargulis}.  
Choose a fine triangulation of  $C_*^{cut}$  
by 3-simplices whose faces 
either are totally geodesic or lie on  
$\partial C_*^{cut}$  so that 
$\varSigma_*^{cut}$  is contained in 
the 1-skeleton.  
Then since  $C_*$  is a geometric limit, 
we may choose a map from the $0$-skeleton of 
$C_*^{cut}$  to   
$C_k$  for sufficiently large  $k$,  where 
4 vertices spanning an oriented simplex in  
$C_*^{cut}$  are mapped 
to 4-vertices spanning a simplex with 
the same orientation also in  
$C_k$  and 
vertices on  $\varSigma_*^{cut}$  are mapped to points
on  $\varSigma_k$.  
Then its obvious piecewise linear extension is necessarily 
an into approximate homeomorphism  
$\varphi_k : C_*^{cut} \to C_k$  
for sufficiently large  $k$. 
\end{pf} 

\begin{LemNum}
$\varphi_k$  can be modified by an isotopy to 
a homeomorphism, which we again denote by  $\varphi_k$,  
so that each component of  $\varphi_k(\partial C_*^{cut})$  
bounds an equidistant tubular 
neighborhood of either  a short geodesic or 
a component of  
$\varSigma_k$  in  $C_k$  
for further sufficiently large  $k$. 
Moreover the isotopy can be chosen so that 
the injectivity radii of the components of 
$\varphi_k(\partial C_*^{cut})$  is uniformly 
bounded from below by some positive constant.    
\end{LemNum} 

\begin{pf} 
Choose a component  $\partial_0 C_*^{cut}$  of  $\partial C_*^{cut}$  
and let  $H_k$  be the image of  $\partial_0 C_*^{cut}$  by  $\varphi_k$,  
namely  $H_k = \varphi_k(\partial_0 C_*^{cut})$.  
Since the nonsingular part  $N_k = C_k - \varSigma_k$  is 
irreducible and atoroidal as a 3-manifold by 
Theorem~\ref{ComplementHyperbolic}, 
$H_k$  either is incompressible and boundary parallel or 
bounds a solid torus  $Z_k$    
in  $N_k$.   
In particular  $H_k$  separates  $N_k$.  

If  $H_k$  is incompressible,  
it is isotopic to a component  
of  $\partial {\cal T}_k$. 
Hence it is isotopic to a horotorus bounding 
an equidistant tubular neighborhood of  
a component of  $\varSigma_k$. 
Choose a horotorus  ${\cal H}_k$  isotopic to 
a corresponding component of  $\partial {\cal T}_k$  
so that 
\begin{equation*} 
\inj {\cal H}_k = I_3/2.   
\end{equation*} 
It exists certainly in  ${\cal T}_k$  
since the minimum of the injectivity 
radius of components of  $\partial {\cal T}_k$  is  
$\geq \inf_{f_k} C_k \geq I_1 \geq I_3$.   
Moreover, since  $\varphi_k$  does not change 
injectivity radius very much, 
it is contained outside the image of  
$C_*^{cut}$  by  $\varphi_k$.  
Hence we can choose isotopy of  $H_k$  to  
${\cal H}_k$  by pushing  $H_k$  outside  
$\varphi_k(C_*^{cut})$.  
In particular, the isotopy is covered by 
an isotopy of  $C_*^{cut}$  fixing the 
complement of a small collar neighborhood of  $\partial_0 C_*^{cut}$,  
and the covering isotopy does not affect the other component 
of  $\partial C_*^{cut}$.  

Suppose next that  $H_k$  bounds a solid torus  $Z_k$  
in  $N_k$.  
We will show that the solid torus  $Z_k$  bounded 
by  $H_k$  contains a simple closed geodesic 
isotopic to the core of  $Z_k$.  
To see this, 
we will first extend  $\varphi_k$  to an embedding  
$\hat{\varphi}_k$   
of the union of  $C_*^{cut}$  and 
a collar  $F$  of  $\partial_0 C_*^{cut}$  which lies in 
the complement of  $C_*^{cut}$.   
The choice of  $F$  is rather technical and 
will be made below.  
Set 
\begin{equation*} 
I_4^k = \inf \{ \inj_x C_k \, \vert \, x \in H_k \}.  
\end{equation*} 
This constant depends on  $k$  in fact, 
however since injectivity radius does not 
change very much by an approximation  $\varphi_k$, 
we may assume that  $I_4^k$  is bounded from below by 
some positive constant for sufficiently large  $k$.  
Then we set a positive constant  $I_4$  by  
\begin{equation*} 
I_4 = \min \{ \inf \{ I_4^k \}, I_3 \}.  
\end{equation*} 
We now choose a collar  $F$  of  $\partial_0 C_*^{cut}$   so that 
the second shortest geodesic on the boundary component  
$\partial_1 C_*^{cut}$  of $F$  other than  $\partial_0 C_*^{cut}$  
has length  $\leq I_4/2$.  
The same argument of the previous lemma 
shows that the extension  $\hat{\varphi}_k$  
over  $C_*^{cut} \cup_{\partial_0 C_*^{cut}} F$  
exists by taking further sufficiently large  $k$.    

Let  $\ell_1$  and  $\ell_2$  be the shortest two 
geodesics on  $\partial_1 C_*^{cut}$.  
They are not homotopic each other because 
$\partial_1 C_*^{cut}$  has an euclidean structure.  
Since  $\hat{\varphi}_k$  does not change the 
length very much,  
we may assume that the length of  
$\hat{\varphi}_k(\ell_1)$  and 
$\hat{\varphi}_k(\ell_2)$  are  $< I_4$.  
The new boundary component  $\hat{H}_k$  
of  $\hat{\varphi}_k(F)$  other than  $H_k$  
is contained in  $Z_k$.  
Since 
$\hat{\varphi}_k(\ell_1)$  and 
$\hat{\varphi}_k(\ell_2)$  are nonhomotopic 
loops on  $\hat{H}_k$,  
at least one of them, say  
$\hat{\varphi}_k(\ell_1)$,  is homotopically 
nontrivial in  $Z_k$.  
Notice that  $\hat{\varphi}_k(\ell_1)$  has 
length  $< I_4$  and hence 
the length shortening homotopy of  
$\hat{\varphi}_k(\ell_1)$  in 
Lemma~\ref{LengthShorteningHomotopyExist}  
does not go through the point with injectivity 
radius  $\geq I_4/2$.  
On the other hand, any point on  $H_k$  has 
injectivity radius  $\geq I_4^k \geq I_4$.
Thus the length shortening homotopy 
of  $\hat{\varphi}_k(\ell_1)$  stays in  $Z_k$  and   
$\hat{\varphi}_k(\ell_1)$  shrinks to some nonzero multiple of 
a closed geodesic  $\ell$  in  $Z_k$.  
Then by Theorem~\ref{ComplementHyperbolic}, 
$C_k - (\varSigma_k \cup \ell)$  is atoroidal,  
and in particular 
$H_k$  is parallel to a torus bounding
an equidistant tubular neighborhood of  $\ell$.  

We would like to isotope  $H_k$  to an 
equidistant torus to  $\ell$  with uniform 
intrinsic injectivity radius.  
To see this, consider an increasing family of 
equidistant tubular neighborhoods of  $\ell$  in  $Z_k$.    
At some first critical radius,   
the boundary of a neighborhood hits either  $H_k$  or itself.   
If the later is the case, 
choose a homotopically nontrivial loop  $g$  
in the critical 
equidistant tubular neighborhood 
which passes the contact point  $p$.  
$g$  is homotopic to some nonzero multiple 
of  $\ell$,  say  $\ell^d$,  since it lies in the solid 
torus  $Z_k$  and  $p$  is a contact point of 
the critical equidistant tubular neighborhood.  
Then consider the developed image of 
the critical equidistant tubular neighborhood.  
The preimage of  $p$  contains points which cannot be 
joined by the action of  $\rho_k(\ell)$  since 
$p$  is a contact point.  
Simultaneously, 
it must be an orbit of a point 
by the action of  $\rho_k(\ell^d)$  
since  $g$  is homotopic to  $\ell^d$.    
This is impossible.  

Hence the equidistant tubular neighborhood of  $\ell$  
grows up to the one  ${\cal H}_k$  which touches  $H_k$.  
Since each point of  $H_k$  has injectivity radius  
$\geq I_4$  in  $C_4$,  we have 
\begin{equation*} 
\inj {\cal H}_k > I_4. 
\end{equation*} 
Moreover, it is contained outside the image of  
$C_*^{cut}$  by  $\varphi_k$.  
Hence we can choose an isotopy of  $H_k$  to  
${\cal H}_k$  by pushing  $H_k$  into  $Z_k$.  
In particular, the isotopy is covered by 
an isotopy of  $C_*^{cut}$  fixing the 
complement of a small collar neighborhood of  
$\partial_0 C_*^{cut}$.  
The covering isotopy does not affect the other component 
of  $\partial C_*^{cut}$  and we are done. 
\end{pf}

\begin{LemNum}\label{LengthBounded} 
For any geodesic loop  $\ell$  on  $\partial C_*^{cut}$,  
$\length \varphi_k(\ell)$  is bounded by some positive constants 
from both above and below for all  $k$,  where 
$\varphi_k$  is a modified one in the previous lemma.  
\end{LemNum} 

\begin{pf} 
The existence of a lower bound is a 
simple corollary to the previous lemma.  
Recall that the isotopy we constructed 
pushes the boundary into an 
equidistant tubular neighborhood.  
Thus it is distance decreasing on the boundary.  
Hence  
since the original $\varphi_k \vert_{\partial_0 C_*^{cut}} : 
\partial_0 C_*^{cut} \to H_k$  
is an approximation and 
does not change the distance very much 
even when  $k$  varies,  
a uniform upper bound exists.  
\end{pf}

By taking a further subsequence and rearranging
reference numbers  
of the components of  $\varSigma$  if necessary,  
we may assume that 
the $j$th component of  $\partial C_*^{cut}$  is mapped by  
$\varphi_k$  to a torus bounding an equidistant tubular neighborhood 
of  $\varSigma^j_k$  for  
$0 \leq j \leq s$  and of a 
short geodesic 
in  $C_k$  for $s < j \leq t$  
for sufficiently large $k$.  

The following lemma finishes the proof or 
Proposition~\ref{ThickProp2}.  

\begin{LemNum} 
There are no components of  $\partial C_*^{cut}$  
which are mapped by  $\varphi_k$  to  
the torus bounding an equidistant tubular neighborhood 
of a short geodesic.  
In other words, 
$s = t$, and there are no cusp openings 
away from the singularity.  
\end{LemNum} 

\begin{pf} 
Assume contrarily that   $s < t$.  
Filling the $j$th component of the boundary of 
$C_*^{cut}$  
for each  $s < j \leq t$  
by an equidistant tubular neighborhood of a short geodesic 
which the  $j$th component of  
$\varphi_k(\partial C_*^{cut})$  bounds in  $C_k$,  we obtain 
a cone-manifold homeomorphic to  
$C - \cup_{j \leq s} \varSigma^j$.  
In other words, for each sufficiently large  $k$,  
$C_*^{cut}$  produces 
an isometric hyperbolic cone-manifold by some Dehn filling 
on the last  $t-s$ components of  $\partial C_*^{cut}$.   

If the number of slopes appeared in this Dehn filling 
on the $j$th component is finite even as  $k$  varies, 
then by taking a further subsequence, 
we may assume that the slope is unique 
and does not depend on  $k$.  
Denote the geodesic representative of a slope 
on  $\partial C_*^{cut}$  by  $m$.  
Then 
by Lemma~\ref{LengthBounded}, 
$\length \varphi_k(m)$  is bounded from above 
by some constant not depending on  $k$.  
On the other hand, since the 
equidistant tubular neighborhood of 
a short geodesic bounded 
by the $j$th component of  $\varphi_k(\partial C_*^{cut})$  approaches 
a cusp, and since the cone angles of short geodesics stay  $2\pi$,  
the lengths of curves on 
$\varphi_k(\partial C_*^{cut})$  bounding a disk in tubes 
must diverge by Lemma~\ref{CuspOpenning}. 
This is a contradiction.  

Hence the number of slopes appeared in the Dehn fillings  
must be infinite for each component when  $k$  varies. 
Let us reconsider the situation not by cutting toral ends 
but by removing the singularity.  
Come back to the limiting cone-manifold  $C_*$.  
The nonsingular part  $N_* = C_* - \varSigma_*$  
admits a hyperbolic structure by Theorem~\ref{ComplementHyperbolic}.  
We denoted it by  $C_{*,comp}$.  
Then for each  $k$,  
$C_{*,comp}$  produces 
$C_{comp}$   by some Dehn filling 
on the last  $t-s$ cuspidal components of  $C_{*,comp}$.   
The filling slope on each component varies infinitely 
many as  $k$  varies.  
Hence such filling slopes 
accumulate to  $\infty$  which corresponds to the 
complete structure  $C_{*,comp}$.  
This means that $C_{*, comp}$  produces  
$C_{comp}$  by infinitely many hyperbolic Dehn fillings. 
On the other hand, 
the slopes of hyperbolic Dehn fillings
on a fixed hyperbolic manifold which produce the same manifold 
is only finitely many.  
This can be verified for instance by listing 
volumes (see \cite{ThurstonNote}).  
Hence we get a contradiction.  
\end{pf}


\subsection{Algebraic Limits}

We assume the conclusion of  Proposition~\ref{ThickProp1}  
throughout this subsection  
and prove Proposition~\ref{ThickProp3} by lemmas below.  
We continue to use  
$C_*^{cut}$  in the previous subsection.  
$C_*^{cut}$  has 
$s$  toral boundaries and each  
$\varphi_k: C_*^{cut} \to C_k$ 
maps the $j$th component of  $\partial C_*^{cut}$  
to a torus bounding an equidistant tubular neighborhood 
of  $\varSigma_k^j$.  
The cusp opening on  $C_*$  does not simply imply 
$(\beta^j =) \, \lim_{k \to \infty} \theta^j_k = 0$  
for  $j \leq s$  as we pointed out 
in section 3.3. 
However in this case, we have 

\begin{LemNum}\label{Beta=0} 
$\beta^j = 0$  for all $1 \leq j \leq s$.  
\end{LemNum}

\begin{pf}  
Suppose contrarily that  
$\beta^j > 0$  for some  $1 \leq j \leq s$  and 
choose a meridional 
element  $m_k$  bounding  $\varSigma_k^j$ 
on  $\varphi_k(\partial C_*^{cut})$ 
for each  $k$.  
If  
$\{ \varphi_k^{-1}(m_k) \, \vert \, k \}$  contains only finitely many 
isotopy classes of the curves on  $\partial C_*^{cut}$,   
then   $\{ \length \varphi_k^{-1}(m_k) \, \vert \, k \}$  is bounded  
by Lemma~\ref{LengthBounded}.  
On the other hand, since a singular solid torus bounded 
by this component of  $\varphi_k(\partial C_*^{cut})$  approaches 
a cusp, the lengths of meridional elements 
on the boundary must diverge by 
Lemma~\ref{CuspOpenning}.  
This is a contradiction.  

Thus it is enough to show that 
$\{ \varphi_k^{-1}(m_k) \, \vert \, k \}$  contains only finitely many 
isotopy classes of the curves on  $\partial C_*^{cut}$.  
Composing reference homeomorphisms of  \eqref{SettingOne}  in 
\ref{Setting}, we obtain a homeomorphism  
$\xi_{k} \circ \xi_{k_0}^{-1} : C_{k_0} \to C_{k}$. 
Then  $\xi_k \circ \xi_{k_0}^{-1}(m_{k_0})$  is isotopic 
to  $m_k$.  
Fixing  $k_0$  and running  $k > k_0$, 
we have infinitely many homeomorphisms  
$\psi_k = 
\varphi_k^{-1} \circ \xi_k \circ \xi_{k_0}^{-1} \circ \varphi_{k_0} : 
C_*^{cut} \to C_*^{cut}$.  
Since the interior of  $C_*^{cut}$  is 
homeomorphic to  $C_*$,  
$C_*^{cut}$  admits only finitely many isotopy classes 
of automorphisms by Corollary~\ref{AutomorphismFinite}.  
Hence there are only finitely many isotopy classes in 
$\{ \psi_k \, \vert \, k > k_0 \}$.  
This is enough since  
$\psi_k(\varphi_{k_0}^{-1}(m_{k_0}))$  
is isotopic to  $\varphi_k^{-1}(m_k)$.  
\end{pf}

\begin{LemNum} 
The angle set of  $C_*$  is equal to  $B$.  
\end{LemNum} 

\begin{pf}  
Lemma~\ref{Beta=0} shows that  $\beta^j = 0$  for 
$j \leq s$.  
Hence it is equal to the $j$th component of the 
angle set of  $C_*$  since the component corresponds 
to a cusp.  

The other component corresponds to a component of 
$\varSigma_*$  in  $C_*$.  
As we have seen in the proof of Lemma~\ref{NoncompactLimit}, 
a tubular neighborhood of  $\varSigma_*^j$  for  
$j > s$  is a limit of  ${\cal T}_k^j$'s  whose 
cone angles are  $\alpha_k^j$'s.  
Thus the cone angle of  $\varSigma_*^j$  in this 
case is  
$\lim_{k \to \infty} \alpha^j_k$  which equals 
$\beta^j$  by definition. 
\end{pf}

\begin{LemNum}
A sequence  $\{ \rho_k \}$  of holonomy representations of  
$\{ C_k \}$  contains a subsequence converging algebraically 
to the holonomy representation of  $C_*$  with respect 
to the identification by  $\xi_k$.  
\end{LemNum}

\begin{pf}
We eventually obtained an into homeomorphism 
$\varphi_k : (C_*^{cut}, \varSigma_*^{cut}) 
\to  (C_k - \cup_{j \leq s} \varSigma_k^j, \cup_{j > s} \varSigma_k^j)$    
for sufficiently large  $k$  by the lemmas in the 
previous subsection.  
On the other hand, there are reference homeomorphisms 
$\xi_k : C \to C_k$  of  \eqref{SettingOne}  in  \ref{Setting}.  
Then the composition  
$\xi_k^{-1} \circ \varphi_k : 
(C_*^{cut}, \varSigma_*^{cut}) \to 
(C  - \cup_{j \leq s} \varSigma^j, \cup_{j > s} \varSigma_k^j)$   
is an into homeomorphism.  
There are only finitely many isotopy classes of 
such maps since otherwise,  $C$  would admit an 
infinite automorphism group, contradicting 
Corollary~\ref{AutomorphismFinite}.  
Hence 
taking a further subsequence, 
we may assume that  $\xi_k^{-1} \circ \varphi_k$'s   
are isotopic for all  $k$.  
Then the algebraic convergence is a consequence 
of a geometric convergence.  
\end{pf}


\section{Thin Tube} 

In this section, 
we study what happens to a geometric limit  $C_*$  of 
a sequence of deformations   $\{ C_k \}$   in  \ref{Setting} 
when the minimum of maximal tube radii goes to zero.  
The analysis involves noncompact euclidean 3-cone-manifolds with 
noncompact singularity, whose definition would be 
obvious.


\subsection{Brief Outline}

The minimum of  $\radius {\cal T}_k^j$  is 
attained by the $c$th component  ${\cal T}_k^c$  which 
contains the first contact point  $c_k$.  
The underlying assumption throughout this section is 

\begin{AssOwnNum}\label{ThinTube} 
If  $k \to \infty$,  then 
\begin{equation*}  
 \radius {\cal T}_k^c \to 0    
\end{equation*} 
\end{AssOwnNum} 

Under the Assumption~\ref{ThinTube} above,  
we discuss possible degenerations in 
two propositions below. 
The technical assumption there will 
be satisfied by some natural setting which 
we will use later on.   

\begin{PropNum}\label{ThinProp1} 
Under the Assumption~\ref{ThinTube}, 
if there is a constant  $V_{min} > 0$  
such that  $\volume C_k \geq V_{min}$,  and 
if $\beta^j$  is strictly less than  $\pi$  
for all $1 \leq j \leq n$, 
then  $C_*$  is isometric to 
either the euclidean line  $\euc$  
or the half line  $\euc_{\geq 0}$. 
\end{PropNum}

\begin{RemNoNum} 
We do not know whether  $\euc_{\geq 0}$  really occurs 
as a geometric limit.  
\end{RemNoNum} 

\begin{PropNum}\label{ThinProp2} 
If furthermore  $\beta^c > 0$, then  
the rescaling limit  $\overline{C}_*$  of  $\{ C_k \}$  by 
normalizing the radius of  ${\cal T}_k^c$  to be $1$  
is an euclidean cone-manifold isometric to 
$\sphere^2(\alpha, \beta, \gamma) \times \euc$, 
where  $\sphere^2(\alpha, \beta, \gamma)$  is an euclidean 
2-cone-manifold over the 2-sphere 
with three cone points of cone angles  
$\alpha, \beta, \gamma$  such that  
$0 < \alpha, \beta, \gamma < \pi$  and    
$\alpha + \beta + \gamma = 2 \pi$. 
\end{PropNum}

\begin{RemNoNum} 
It is quite unlikely that both  $\radius {\cal T}_k^c \to 0$  and 
$\beta^c = 0$  occur simultaneously, though 
we do not have a proof.  
\end{RemNoNum}


\subsection{Collapsing} 

Consider the Dirichlet polyhedron  
$P_{c_k}$  of  $C_k$  about the first 
contact point  $c_k$,  
which we simply denote by  $P_k$  from now on.  
Supposing that the Assumption~\ref{ThinTube} is the case, 
we analyze the limit of  $P_k$  and 
prove Proposition~\ref{ThinProp1}.  

\begin{LemNum} 
Under the assumption of Proposition~\ref{ThinProp1}, 
$\{ (P_k, c_k) \}$  converges geometrically 
to the euclidean line 
$\euc$  or the half line  $\euc_{\geq 0}$.   
\end{LemNum} 

\begin{pf}
Imagine that  $c_k$  is the contact point of  
${\cal T}_k^c$  from the top side 
and  ${\cal T}_k^{c'}$   from the bottom sides.  
The reference number  $c$  might be equal to  $c'$.  
The shortest common orthogonal to  $\varSigma_k^c$  and  
$\varSigma_k^{c'}$  
which goes through  $c_k$  lifts to 
the geodesic segment  $g_k \subset P_k$.  
It is in fact the segment realizing  
the length $ = 	2 \radius {\cal T}_k^c$, 
and by the Assumption~\ref{ThinTube},  
$\length g_k$  goes to zero as  $k \to \infty$.  

Let  $p_k, q_k$  be terminal points of 
$g_k$.  
Since we assumed that  
$\beta^j$  is strictly less than  $\pi$,  
$P_k$  is locally bounded by roof shaped faces 
near  $p_k$  from the top and  $q_k$  from the bottom 
respectively,  where their ridges correspond to 
$\varSigma_k^c$  and  $\varSigma_k^{c'}$.  
$P_k$  is convex, and it is bounded by 
the extension of 
these roofs from the top and bottom.  
Moreover since  $\length g_k \to 0$, and 
since the volume is assumed to be bounded away from zero, 
the roof ridges become arbitrary close and parallel.  
Hence  $\{ P_k \}$  converges as a metric 
space to a connected closed subset of the euclidean line  $\euc$. 
\end{pf}

\begin{pf*}{\it Proof of Proposition~\ref{ThinProp1}}  
Choose for each  $k$  a segment  $l_k \subset \euc$  
through  $c_k$  so that    
it is maximally embedded in  $P_k$.  
By the previous lemma, we have  
$\length l_k \to \infty$  when 
$k \to \infty$.  
Thus a long segment  $l_k$  
can be isometrically embedded in  $C_k$.  
Assigning to each point of  $C_k$  the 
nearest point on the image of  $l_k$, we 
obtain a map  $\varphi_k : C_k \to l_k \subset \euc$.   
Then the relation  $\text{R}_k$  between  
$C_k$  and  $\euc$  defined by  
$\text{R}_k = \{ (x, y) \in C_k \times \euc \vert 
	\varphi_k(x) = y \}$ 
is an approximation for some  $\varepsilon$  
where  $\varepsilon \to 0$  as  $k \to \infty$.  
\end{pf*}


\subsection{Rescaling}  

In this subsection, under the conclusion of 
Proposition~\ref{ThinProp1}, 
we prove Proposition~\ref{ThinProp2}.  

\begin{LemNum}\label{RescalingLimit}
If  $\beta^c > 0$,  then 
the rescaling limit of  $\{ C_k \}$  normalizing 
the radius of  ${\cal T}_k^c$  to be  $1$  is a noncompact 
euclidean cone-manifold with 
nonempty singular set.  
\end{LemNum}

\begin{pf} 
Multiply  
$1/\radius {\cal T}_k^c$  on the metric of  $C_k$,  
we obtain a cone-manifold  
$\overline{C}_k$  of constant 
curvature  $= -(\radius {\cal T}_k^c)^2$,  
which is  $\geq -1$  for large  $k$.  
Then 
$\radius \overline{{\cal T}}_k^c = 1$. 
Notice that the estimate \eqref{InjBound} is in fact 
valid for tubes with constant sectional curvature $= K$  
where  $-1 \leq K \leq 0$,  because the bound in the  
hyperbolic case 
is the worst.  
Since  $\beta^c$  is strictly less than  $\pi$,  
$\inj \partial \overline{{\cal T}}_k^c \leq \pi \sinh 1$  
by this new estimate.  
Then by Lemma~\ref{InjectivityComparison}, we have a constant 
$W = W(1, \pi \sinh 1, -1)$  such that 
\begin{equation*} 
\inj \partial \overline{{\cal T}}_k^c \leq 
	W \cdot \inj_{\overline{c}_k} \overline{C}_k.  
\end{equation*} 
Since the $c$th component  $\beta^c$  of the 
angle set is positive by the assumption, and 
$\length \overline{\varSigma}_k^c$  
diverges, 
we can embed an euclidean disk of radius  
$\beta^c/2$  into  $\partial \overline{{\cal T}}_k^c$  by 
the euclidean case of \eqref{Rectangle}.  
Therefore,     
$\inj \partial \overline{{\cal T}}_k^c \geq \beta^c/2$  
for sufficiently large  $k$,     
$\inj_{\overline{c}_k} \overline{C}_k$  is 
uniformly bounded from below,  and 
in particular, $\overline{c}_*$  admits an euclidean 
ball neighborhood. 

On the other hand,  
each nonsingular component  $\varSigma_k^j$  either becomes 
parallel to  $\varSigma_k^c$  in  $C_k$  or  goes far away 
from   $c_k$.  
In particular, either  $\length \varSigma_k^j \to \infty$  
or  $d(c_k, \varSigma_k^j) \to \infty$.  
This is true also in the rescaled setting.  

These two informations are good enough to conclude that  
$\{ (\overline{C}_k, \overline{c}_k) \}$  converges 
geometrically to an euclidean 
cone-manifold  $(\overline{C}_*, \overline{c}_*)$  since 
the singularity 
admits uniformly thick tubular neighborhood and 
its length does not degenerate, also the reference 
point stays in a uniformly thick part.  
The singular set is nonempty because  
$\overline{\varSigma}_*^{c}$  has distance  $1$  to  
$\overline{c}_*$.  
\end{pf}

\begin{LemNum} 
$\overline{C}_*$  has two ends.
\end{LemNum} 

\begin{pf} 
That  
$\overline{C}_*$  has two ends is equivalent to that 
$C_*$  converges to  $\euc$  instead 
of  $\euc_{\geq 0}$.  

Assume contrarily that  $C_*$  converges to  
$\euc_{\geq 0}$  and 
choose  $R > 0$  sufficiently large so that  
$\partial \ball_R(\overline{C}_*, \overline{c}_*)$  is connected.  
$\partial \ball_R(\overline{C}_*, \overline{c}_*)$  can be seen 
in the Dirichlet polyhedron  
$\overline{P}_*$  of  $\overline{C}_*$  about  $\overline{c}_*$  
as an intersection of  $\overline{P}_*$  
and the sphere of radius  $= R$.  
The faces of  $\overline{P}_*$  
intersecting  
$\partial \ball_R(\euc^3, \overline{c}_*)$  for large  $R$  all 
must be parallel to the ray to the end,  
and hence the combinatorial structure of  
$\partial \ball_R(\overline{C}_*, \overline{c}_*)$  is very simple. 
It consists of one 2-cell  $e$  with  $\mu$  edges,  where 
$\mu$  is equal to the number of faces of  $\overline{P}_*$  
intersecting  $\partial \ball_R(\euc^3, \overline{c}_*)$.  
Also the topology of  
$\partial \ball_R(\overline{C}_*, \overline{c}_*)$  does not change 
for sufficiently large  $R$  
since  $\partial \ball_R(\overline{C}_*, \overline{c}_*)$  tends 
to be orthogonal to the ray to the end. 

To see more about vertices, we let  $\nu_i$  be 
the total angle of corners of  $e$  surrounding the $i$th 
vertex of  $\partial \ball_R(\overline{C}_*, \overline{c}_*)$  
and  $\nu$  the number of vertices.  
Then by Gauss-Bonnet, we have the identity.  
\begin{align*} 
- \int_{\partial \ball_R(\overline{C}_*, \overline{c}_*)} K_R \, dA 
	& = \int_{\partial e} \kappa_g \, ds  
		+ (\mu-2) \pi - \sum_{i=1}^{\nu} \nu_i \\
	& = \int_{\partial e} \kappa_g \, ds 
		- 2 \pi \chi(\partial \ball_R(\overline{C}_*, \overline{c}_*)) 
			+ \sum_{i=1}^{\nu} (2\pi - \nu_i), 
\end{align*}
where  $K_R$  is a Gaussian curvature of  
$\partial \ball_R(\overline{C}_*, \overline{c}_*)$  supported 
on the interior of  $e$  and  $\kappa_g$  
is a geodesic curvature along  $\partial e$.  

Let us see what happens when  $R \to \infty$.  
The left hand side goes to zero since  $K_R \to 0$  and 
$\area \partial \ball_R(\overline{C}_*, \overline{c}_*)$  is bounded.  
The first term of the right hand side approaches zero also since 
$\kappa_g \to 0$  but the length of  $\partial e$  is 
bounded.  
For each vertex not on  $\varSigma_*$,  $\nu_i \to 2\pi$,  
and on  $\varSigma_*$, $\nu_i \to \beta^j$  where 
$\beta^j$  is a cone angle 
of the singularity on which 
the $i$th vertex lies.  
Hence if  $R$  is large enough, the contribution of  
$2 \pi \chi(\partial \ball_R(\overline{C}_*, \overline{c}_*))$,  which 
is constant, and the contribution of the cone angles of  
$\varSigma_* \cap \partial \ball_R(\overline{C}_*, \overline{c}_*)$  are identical.  

Since we have assumed that  $0 \leq \beta^j < \pi$  for all  $j$  
but  $c$  and  $0 < \beta^c < \pi$, 
this cancellation occurs only when  
$\chi(\partial \ball_R(\overline{C}_*, \overline{c}_*)) = 2$  and  
$\varSigma_* \cap \partial \ball_R(\overline{C}_*, \overline{c}_*)$  
consists of three points.  
This is a contradiction since a noncompact euclidean 
cone-manifold with one end must have 
even number of ends of singularity.  
\end{pf} 

The following classification of 
noncompact euclidean 3-cone-manifolds   
with two ends whose cone angles all are  $\leq \pi$  
finishes the proof of Proposition~\ref{ThinProp2}.  

\begin{LemNum}\label{Classification}
An orientable euclidean 3-cone-manifold  $E$  with 
nonempty singular set of cone angles 
$\leq \pi$  and with two ends is 
a product of a compact euclidean 2-cone-manifold 
and  $\euc$.  
More precisely, $E$  is isometric to either 
\begin{enumerate} 
\item $\sphere^2(\pi, \pi, \pi, \pi) \times \euc$ \; or
\item $\sphere^2(\alpha, \beta, \gamma) \times \euc$,  
	where  $\alpha + \beta + \gamma = 2 \pi$.  
\end{enumerate} 
\end{LemNum}

\begin{pf}
Choose a reference point  $e \in E$.  
A Dirichlet polyhedron  $P_e$ 
is a convex polyhedron.  
Since  $E$ has two ends, there are two rays  
$r_1, r_2$  in  $P_e$  based at  $e$.  
If  $r_1 \cup r_2$  had bent at  $e$,  then 
$P_e$  cannot have two ends.  
Hence  $r_1 \cup r_2$  is a straight line.  
Then any face  $P_e$  must be parallel to 
$r_1 \cup r_2$  by convexity and they surround  
$r_1 \cup r_2$.  
Let  $Q$  be a polygon through  $e$  which intersects 
perpendicularly to the faces.  
It must be compact since otherwise, 
$P_e$  would not have two ends.  
Then  $\partial Q$  is glued 
with  $\partial Q$ via the identification of 
$P_e$  because  $P_e$  is starlike and 
hence the identifications do not contain any 
translation factor along  $r_1\cup r_2$.  
Thus  $Q$  becomes an euclidean sub-cone-manifold 
after identification and 
$E$  is the product of this sub-cone-manifold 
and  $\euc$.  

The remaining is to classify compact 
euclidean 2-cone-manifolds with 
cone angles  $\leq \pi$.  
However this is a routine application of 
Gauss-Bonnet, so that 
we leave it to the reader.   
\end{pf}


\section{Continuous Families} 

In this section, we come back to a continuous family of 
deformations of a compact orientable 
hyperbolic 3-cone-manifold  $C$,  see what 
happens in the limit 
with the aid of the propositions in the 
previous sections, and prove the main theorem 
and its corollaries.


\subsection{Brief Outline}

Let  $\{ C_{\theta} \}$  be 
a continuous family of deformations of  $C$  
parameterized by the angle assignment  
\begin{equation*}
\theta : [0, 1) \to \Theta = (\Im \complex)^n
\end{equation*} 
where 
\begin{equation*} 
\lim_{t \to 1} \theta(t) = B = (\beta^1, \cdots, \beta^n).  
\end{equation*} 
We will first of all 
generalize the concept of strong convergence 
of a sequence to a continuous family of deformations.  

\begin{DefNoNum} 
A continuous family  $\{ C_{\theta} \}$  of the 
deformations of  $C$  is said to 
converge strongly to  $C_*$  if every subsequence  
$\{ C_k \}$  in  $\{ C_{\theta} \}$  whose angle set 
tends to  $B$  converges strongly to  $C_*$.  
\end{DefNoNum} 

The following theorems, which we will prove in this section, 
are what we can conclude for
angle monotone families from the analysis of the 
previous sections.  
As we will see later in the proof, 
the strong convergence of a family is very much likely derived 
from a strong convergence of a sequence together 
with the local rigidity.  

\begin{ThNum}\label{AngleDecreasingDeformation} 
If the family  $\{ C_{\theta} \}$  has a 
component-wise decreasing angle assignment  $\theta$    
and if  $\beta^j > 0$  for all  $1 \leq j \leq n$,  
then $\{ C_{\theta}\}$  converges strongly to 
a compact hyperbolic cone-manifold  $C_*$  
homeomorphic to  $C$.  
\end{ThNum}  

\begin{RemNoNum} 
The positivity of  $\beta^j$   is conjecturally 
unnecessary.  
However, the proof we present depends on Proposition~\ref{ThinProp2} 
which involves this unclear hypothesis.  
Also it forces us to make a technical arrangement 
in the proof of the main theorem.  
\end{RemNoNum} 

\begin{ThNum}\label{AngleIncreasingDeformation}
If the family  $\{ C_{\theta} \}$  has a 
component-wise increasing angle assignment  $\theta$    
and  $\beta^j < \pi$  for all  $1 \leq j \leq n$,  
then either 
\begin{enumerate} 
\item\label{AIDOne} 
	$\volume C_{\theta} \to 0$, 
\item\label{AIDTwo}  
	$C_{\theta}$  contains 
	a sub-cone-manifold homeomorphic to  $\sphere^2$  with three 
	cone points so that the sum of their cone angles 
	approaches  $2\pi$  as  $\theta \to B$,  or   
\item\label{AIDThree} 
	$\{ C_{\theta} \}$  converges strongly to 
	a hyperbolic cone-manifold homeomorphic to  $C$.  
\end{enumerate} 
\end{ThNum}
  
The first theorem will be used to prove the 
main results in the final subsection.  
The second theorem above is just to note what 
we can conclude for angle increasing family, 
and it is not related to the main results directly.


\subsection{Splitting} 

The following example shows the 
degeneration discussed in  Proposition 
\ref{ThinProp2}  really occurs in an 
angle increasing continuous family.  

\begin{ExamNoNum} 
Let  $\alpha, \beta, \gamma$  be numbers such that  
$0 < \alpha, \beta, \gamma < \pi$  and  
$\alpha + \beta + \gamma = 2 \pi$.  
Consider a hyperbolic tetrahedron whose 
dihedral angles along ridges are  
$(\alpha - \varepsilon)/2, \, 
(\beta - \varepsilon)/2, \, (\gamma - \varepsilon)/2$  in 
three opposite pairs,  
where  $\varepsilon$  is a small nonnegative number.  
When  $\varepsilon = 0$,  the tetrahedron  
has four ideal vertices.  
If  $\varepsilon > 0$,  then the tetrahedron is of infinite 
volume.  
It becomes compact by 
truncating the ends by polar planes.  
The result is called a truncated tetrahedron.  
We thus obtain a family of polyhedra   
$\{ \Delta_{\varepsilon}(\alpha, \beta, \gamma) \}$  of 
finite volume 
parameterized by  $\varepsilon \geq 0$.  
    
When  $\varepsilon > 0$,  
taking a double of  
$\Delta_{\varepsilon}(\alpha, \beta, \gamma)$  along 
$4$ hexagonal faces, we obtain a compact 
hyperbolic cone-manifold with geodesic boundary.  
Taking further double along the boundary, 
we obtain a family of closed hyperbolic 
3-cone-manifolds  $\{ C_{\varepsilon}(\alpha, \beta, \gamma) \}$.  
The singular set  $\varSigma$  consists of 
6 circles each of which is assigned 
$\alpha - \varepsilon, \, 
\beta - \varepsilon, \, \gamma - \varepsilon$  in pairs 
as their cone angles.  
If  $\varepsilon \to 0$, then 
the face of truncation becomes very tiny, and  
$C_{\varepsilon}(\alpha, \beta, \gamma)$  splits into two 
noncompact hyperbolic cone-manifolds 
by tearing off the boundary 
of the first double.  
This family shows that 
the degeneration in Proposition~\ref{ThinProp2}  
certainly occurs at the face of truncation.  
The reference point lies on the 
face of truncation, 
and the rescaling geometric limit 
is isometric to  
$\sphere^2(\alpha, \beta, \gamma) \times \euc$.  
\end{ExamNoNum}

This example gives us a fairly general picture.   
In fact, using the following observation,  
we will show in the proof of  Theorem~\ref{AngleDecreasingDeformation} 
that the splitting degeneration caused by 
appearance of an euclidean sub-cone-manifold such as this 
cannot occur in angle decreasing families.  

\begin{LemNum}\label{AxisMeet} 
Let  $\iota_1, \iota_2$  be elliptic elements 
in  $\PSL_2(\complex)$  with axis  
$\ell_1, \ell_2$  respectively.  
If  $\iota_1 \iota_2$  is elliptic with axis  $\ell_{12}$, 
and 
if the total angle of rotations of   
$\iota_1, \iota_2$  and  $\iota_1 \iota_2$  is   
$> 2 \pi$,  
then  $\ell_1, \ell_2$  
and  $\ell_{12}$  meets at the unique point in  $\hyp^3$.  
\end{LemNum} 

\begin{pf} 
It is not quite hard to show that 
$\iota_1 \iota_2$  is loxodromic if  
the union  $\ell_1 \cup \ell_2$  does not lie on a geodesic 
plane in  $\hyp^3$.  
Hence we may assume that  $\ell_1 \cup \ell_2$  lies on 
a geodesic plane  $X$.  
Replacing the role of  $\iota_1, \iota_2$  
by  $\iota_1^{-1}, \iota_1 \iota_2$  and  
$\iota_1 \iota_2, \iota_2^{-1}$,  
we get geodesic planes  $Y$  and  $Z$  supporting  
$\ell_1 \cup \ell_{12}$  and  $\ell_{12} \cup \ell_2$ 
respectively.   
If three planes  $X, Y, Z$  meets in  $\hyp^3$, 
then we are done.  
If not, they either meet at  $\sphere^2_{\infty}$  
or does not meet in  $\hyp^3 \cup \sphere^2_{\infty}$  
and admits a geodesic plane meeting  $X, Y, Z$  
perpendicularly.  
In both cases, the sum of three rotation angles 
must be  $\leq 2 \pi$  and 
the assumption is not satisfied.    
\end{pf}


\subsection{Angle Decreasing Family} 

In this subsection, we prove 
Theorem~\ref{AngleDecreasingDeformation} using 
propositions in the previous sections and 
the observation in Lemma~\ref{AxisMeet}.  
First of all, we have a lower bound of the volume.  

\begin{LemNum}\label{Hodgson} 
If $\theta$  is component-wise 
decreasing, then there is 
a constant  $V_{min} > 0$  such that 
\begin{equation*} 
V_{min} \leq \volume C_{\theta(t)} 
\end{equation*} 
for all  $t \in [0, 1)$.  
\end{LemNum}

\begin{pf}
Recall Schll\"affli's variation formula revisited by 
Hodgson \cite{Hodgson}, 
\begin{equation*}\label{??}
d \volume C_{\theta} = 
		- \frac{1}{2} \, \sum_{j=1}^n \, \length \varSigma^j \, 
		d \, \theta^j,   
\end{equation*} 
where  $\theta^j$  is the $j$th component of  $\theta$.  
It says that the volume is an increasing function 
in angle decreasing deformations.  
Hence  
$\volume C_{\theta(t)} \geq \volume C_{\theta (0)} =\volume C$.  
\end{pf}

\begin{pf*}{Proof of Theorem~\ref{AngleDecreasingDeformation}} 
Given an angle decreasing family   $\{ C_{\theta} \}$,  
where  $\beta^j > 0$  for  $1 \leq j \leq n$,  
we set  $C_i = C_{\theta(1-1/i)}$  and choose a geometrically 
convergent subsequence  $\{ C_k \}$  
in \ref{Setting}  with canonical reference 
homeomorphisms  $\{ \xi_k \}$.  

Assume that  $\{ C_k \}$  satisfies the 
Assumption~\ref{ThinTube}.  
Since the volume is bounded from below, and also 
since  $0 < \beta^c < \pi$,  the rescaling geometric limit 
$\overline{C}_*$  is by Proposition~\ref{ThinProp2} 
isometric to the product  
$\sphere^2(\alpha, \beta, \gamma) \times \euc$  
where  $\alpha + \beta + \gamma = 2\pi$.  
Since  $\overline{C}_*$  contains an euclidean 2-cone-manifold 
as a section, we can find a topologically same section by
an approximation in a reasonably large neighborhood of  $c_k$,  
say  $\ball_R(C_k, c_k)$,  for sufficiently large  $k$.  
It is homeomorphic to the 2-sphere transversely intersecting 
$\varSigma_k$  at three points.  
The total sum of cone angles of these points is more 
than  $2\pi$  since the deformation is angle decreasing. 

Now, there are two components  $\ell_1$  and  $\ell_2$  of  
$\ball_R(C_k, c_k) \cap \varSigma_k$  which admit the shortest 
common orthogonal  $g$  going through  $c_k$.  
Develop  $\ell_1 \cup g \cup \ell_2$,  then the images of 
$\ell_1$  and  $\ell_2$  cannot have a common point 
even in their extensions. 
On the other hand, if we let two meridional 
elements rounding  $\ell_1$  and  $\ell_2$  be 
$m_1$  and  $m_2$  respectively, 
then since $m_1 m_2$  becomes a meridional element 
rounding the last component, 
$\rho_k(m_1 m_2)$  represents an elliptic element.  
Moreover the total angles of rotations of  
$\rho_k(m_1), \rho_k(m_2)$  and  $\rho_k(m_1m_2)$  
is  $> 2 \pi$.  
Thus by Lemma~\ref{AxisMeet}, the developed 
image of  $\ell_1$  and  $\ell_2$  must have 
common point in their extensions.  
This is a contradiction.  

Hence  $\{ C_k \}$  does not satisfy 
the Assumption~\ref{ThinTube},  and 
the radius of the maximal tube must be uniformly 
bounded away from zero.  
We can now apply the results in section 5.  
In particular, 
the geometric limit $C_*$  is 
a strong limit of a sequence  $\{ C_k \}$  by 
Corollary~\ref{StrongLimit}.  

To see that 
$C_*$   is a strong limit of a family  $\{ C_{\theta} \}$, 
let  $\rho_*$  be a holonomy representation of  $C_*$.    
Since it is a holonomy representation of 
a cone-manifold  $C_*$,  it can be deformed 
in a small range by Corollary 
\ref{SmallDeformationExist}.  
Let us choose a small path 
on the space of representations   
$\Hom (\Pi, \PSL_2(\complex))/\PSL_2(\complex)$  
from  $\rho_*$  supported on  $[0, \varepsilon)$  
so that the associated angle assignment
is equal to  $\theta(1 - t)$  where  
$t \in [0, \varepsilon)$.  
This path and the path defined by  
$\{ \rho_{\theta(t)} \}_{0 \leq t < 1}$  in the 
space of representations have 
common points accumulating  $\rho_*$,  which are 
realized by holonomy representations  $\{ \rho_k \}$  of  
$\{ C_k \}$  in  \ref{Setting}.    
Then since the paths are the image of the same  
angle assignment, they must be the same by the 
local rigidity at  $\rho_*$.  
\end{pf*}


\subsection{Angle Increasing Family} 

In this subsection, we prove 
Theorem~\ref{AngleIncreasingDeformation}  
and present one example for which the theorem 
can be applied.  

\begin{pf*}{Proof of Theorem~\ref{AngleIncreasingDeformation}}  
Assume that  $\volume C_{\theta}$  does not converge to zero,  
in other words, \eqref{AIDOne}  is not the case.  
Choose a sequence  $\{ (C_k, c_k) \}$  as in  \ref{Setting}.  
If the Assumption~\ref{ThickTube} is the case, 
then by Corollary~\ref{StrongLimit}, we have 
a strong limit  $C_*$  of  $\{ C_k \}$.  
$C_*$  is also a strong limit of a family  $\{ C_{\theta} \}$  
by the same argument in Theorem~\ref{AngleDecreasingDeformation}, 
and we get the third case.  
If the Assumption~\ref{ThinTube} is the case, 
then by Proposition~\ref{ThinProp2}, 
the rescaling limit  $\overline{C}_*$  is isometric 
to  $\sphere^2(\alpha, \beta, \gamma) \times \euc$.  
Thus  $\overline{C}_*$  contains an euclidean 
sub-cone-manifold as its section.  
We then have a topologically same section in 
$C_k$  for sufficiently large  $k$  
by an approximation and hence 
we are in the second case.  
\end{pf*} 

\begin{ExamNoNum} 
This observation can be used for example to 
study an angle increasing family 
$\{ \text{\bf 8}_{\theta} \}$  on the 3-sphere singular 
along the figure eight knot (see \cite{ThurstonNote}).  
Since an underlying space of  $\text{\bf 8}_{\theta}$  
is the 3-sphere, \eqref{AIDTwo} does not occur.  
Hence the angle increasing deformation is possible 
as long as  
$\volume \text{\bf 8}_{\theta} > 0$.  
The $A$-polynomial in \cite{CCGLS},  
which is  
\begin{equation}\label{A-polynomial}  
-M^4 +L(M^8-M^6-2M^4-M^2+1)-L^2M^4    
\end{equation} 
for the figure eight knot for instance, 
represents 
a relation between an eigenvalue  $M$  of an meridian 
and an eigenvalue  $L$  of  
a longitude for  $\SL_2(\complex)$-representations 
of a knot group.  
Then setting  $M = \exp(t \sqrt{-1}/2)$  
in the equation  \eqref{A-polynomial} $= 0$,  
we obtain  
\begin{equation*} 
\cosh \log (-L) = -\frac{L+L^{-1}}{2} = 1 + \cos t -\cos 2 t.  
\end{equation*} 
This shows that  $L$  is always real and the length of the 
singularity at  $t$  is equal to  $2 \log (-L)$.  
Thus by Schl\"affli's formula, 
\begin{align*}  
\volume \text{\bf 8}_{\theta} 
	& = - \int_0^{\theta} \log (-L) d t + \volume \text{\bf 8}_0 \\ 
	& = - \int_0^{\theta} \arccosh(1+\cos t - \cos 2 t) \, d t 
		-6 \int_0^{\pi/3} \log \vert 2 \sin t \vert \, d t 
\end{align*} 
Hence to find the deformable range is reduced to 
the computation of this integral.  
A numerical computation shows that 
$\text{\bf 8}_{\theta}$  
survives as long as  $\theta < 2\pi/3$.  
\end{ExamNoNum}


\subsection{Proof of Theorem and Corollaries} 

\begin{pf*}{\it Proof of Theorem}
Given a compact orientable hyperbolic 
3-cone-manifold  $(C, \varSigma)$  
with an angle set   $A = (\alpha^1, \cdots, \alpha^n)$,   
we start with a complete structure  $C_{comp}$  
supported on the nonsingular part  $N = C - \varSigma$.  
Since a small angle changing deformation of  $C_{comp}$  
uniquely exists by 
Corollary~\ref{SmallDeformationExist}, 
there is an angle set  
$B = (\beta^1, \cdots, \beta^n)$  very 
close to  $(0, \cdots, 0)$  so that 
$C_{comp}$  admits angle increasing 
deformations along 
a linear path  $\zeta : [0,1] \to \Theta$  
with  $\zeta(1) = B$.    
Moreover we can choose 
each  $\beta^j$  positive so that  
$\beta^j = 2\pi/b_j \; (< \alpha^j)$  for some 
large integer  $b_j$.  
Let  $\{ C_{\zeta(t)} \}$  be the associated 
family of deformations.  
The cone-manifold  $C_{\zeta(1)}$  with 
the angle set   $B$  
shares the topological type with  $C_{comp}$  and 
hence  $C$.  
Moreover  $C_{\zeta(1)}$  is an orbifold.    

Choose a linear path  
$\theta : [0,1] \to \Theta$  between 
$A$  and  $B$  such that  
$\theta(0) = A$  and  $\theta(1) = B$. 
It is component-wise decreasing since  
$\beta^j < \alpha^j$.  
We have not known that the path  $\theta$  is 
supported by a continuous family of 
deformations in the full range.  
But since there always exists a small deformation  
by Corollary~\ref{SmallDeformationExist}, 
we may assume that  $C$  is actually deformable 
at least in the range  $[0,\omega)$  for some   $0 < \omega \leq 1$.  
Then since  $\theta$  is angle decreasing, 
the family converges strongly to a compact 
hyperbolic 3-cone-manifold  $C_*$  by 
Theorem~\ref{AngleDecreasingDeformation},  where 
the angle set of  $C_*$  is equal to  
$\theta(\omega)$.  
Then by Corollary~\ref{SmallDeformationExist}, 
the deformation can extend further.    
The prolongation of the deformable range can be done 
up to when  $t$  reaches to  $1$.  
Hence we have obtained a continuous family of 
deformations  
$\{ C_{\theta} \}$  
for full range of  $\theta$.  
$C_{\theta(1)}$  is homeomorphic to  $C$.  
Moreover  $C_{\theta(1)}$  is an orbifold, and 
hence  $C_{\theta(1)}$  and  $C_{\zeta(1)}$  are isometric 
by Mostow rigidity.  
We thus have connected two 
cone-manifolds  $C$  and  $C_{comp}$  
through  $\{ C_{\theta} \}$  and  $\{ C_{\zeta}\}$.  
\end{pf*}

\begin{pf*}{\it Proof of Corollary 1} 
Suppose we are given two cone-manifolds  $C$  and  $C'$  
which are isomorphic.  
They can be deformed along the same 
angle decreasing path used in the proof 
of Theorem to the complete 
manifold  $C_{comp}$  and  $C'_{comp}$.  
The destinations are isometric 
by Mostow-Prasad rigidity.  
Then 
the returning path to  $C$  and  $C'$  
must be the same since the 
cone angle is the only parameter 
by the local rigidity.  
\end{pf*}

\begin{pf*}{\it Proof of Corollary 2}
This is now obvious since 
our family is supported by 
a path of holonomy representations 
and one terminal corresponds to the complete 
structure which is liftable. 
\end{pf*}


\bibliographystyle{amsplain}

\end{document}